\documentclass[11pt]{amsart}

\usepackage{overpic, amssymb, times, graphicx, tikz-cd,float}
\usepackage[all]{xy}
\usepackage[backref=page]{hyperref}
\usepackage{verbatim}
\usepackage{xcolor}
\usepackage{setspace}

\DeclareUnicodeCharacter{2011}{-}
\usepackage{cleveref}
\usepackage{bm}
\usepackage{microtype}
\usepackage[T1]{fontenc}
\usepackage{lmodern}
\usepackage{amsmath,amsthm}
\usepackage{thmtools}
\usepackage{enumitem}

\hypersetup{
    colorlinks=true,
    linkcolor=blue,
    filecolor=blue,      
    urlcolor=blue,
    citecolor=blue,
}
\numberwithin{equation}{section}

\newtheorem{theorem}{\textbf{Theorem}}[section]
\newtheorem*{theorem*}{\textbf{Theorem}}
\newtheorem{thmx}{Theorem}
 
\newtheorem{corollaryx}[thmx]{\textbf{Corollary}}

\newtheorem{definition}[theorem]{\textbf{Definition}}
\newtheorem{proposition}[theorem]{\textbf{Proposition}}
\newtheorem{lemma}[theorem]{\textbf{Lemma}}

\newtheorem{claim}[theorem]{\textbf{Claim}}
\newtheorem{corollary}[theorem]{\textbf{Corollary}}
\newtheorem{remark}[theorem]{\textbf{Remark}}
\newtheorem{assumption}[theorem]{\textbf{Assumption}}
\newtheorem{example}[theorem]{\textbf{Example}}

\newtheorem{definition/proposition}[theorem]{\textbf{Definition/Proposition}}

\AddToHook{env/example/begin}{\pushQED{\qed}}
\AddToHook{env/example/end}{\popQED}

\def\N{{\mathbb N}}
\def\R{\mathbb{R}}
\def\Z{{\mathbb Z}}

\def\C{{\mathbb C}}
\def\D{{\mathbb D}}

\def\Q{{\mathbb Q}}

\def\K{{\mathbb K}}

\def\bbS{{\mathbb S}}
\newcommand{\CP}{\mathbb{C}\mathbb{P}}
\newcommand{\RP}{\mathbb{R}\mathbb{P}}

\def\cA{{\mathcal A}}

\def\cC{{\mathcal C}}

\def\cH{{\mathcal H}}

\def\cL{{\mathcal L}}
\def\cM{{\mathcal M}}

\def\cP{{\mathcal P}}

\def\rD{{\rm D}}

\def\rd{{\rm d}}

\def\la{\langle\,}
\def\ra{\,\rangle}

\def\std{\rm std}

\def\CZ{\rm CZ}

\def\ev{\mathrm{ev}}

\DeclareMathOperator{\Id}{id}

\DeclareMathOperator{\coker}{coker}

\DeclareMathOperator{\rank}{rank}

\DeclareMathOperator{\Tor}{Tor}

\DeclareMathOperator{\CH}{CH}
\DeclareMathOperator{\CC}{CC}

\DeclareMathOperator{\LCH}{LCH}
\DeclareMathOperator{\GW}{GW}
\DeclareMathOperator{\QH}{QH}
\DeclareMathOperator{\RFH}{RFH}
\DeclareMathOperator{\FS}{FS}
\DeclareMathOperator{\pt}{pt}
\DeclareMathOperator{\Proj}{Proj}
\newcommand{\SH}{\operatorname{SH}}
\newcommand{\FH}{\operatorname{FH}}
\newcommand{\FC}{\operatorname{FC}}
\newcommand{\Spec}{\operatorname{Spec}}
\newcommand{\bd}{\partial}

\newcommand{\PD}{\operatorname{PD}}
\newcommand{\Addresses}{{
		\bigskip
		\footnotesize

	    Zhengyi Zhou, \par\nopagebreak
        \textsc{State Key Laboratory of Mathematical Sciences, Chinese Academy of Sciences;}\par\nopagebreak
	    \textsc{Morningside Center of Mathematics, Chinese Academy of Sciences;}\par\nopagebreak
         \textsc{Academy of Mathematics and Systems Science, Chinese Academy of Sciences, China}\par\nopagebreak
		\textit{E-mail address}: \href{mailto:zhyzhou@amss.ac.cn}{zhyzhou@amss.ac.cn}

}}

\title{Symplectic $\CP^n$ recognition and filling rigidity \\ via Rabinowitz Floer homology}
\author{Zhengyi Zhou}

\begin{document}
	\maketitle
\begin{abstract}
We prove that a symplectic manifold $X$ admitting a pseudo-holomorphic line through two generic points, formulated as the non-vanishing of a genus-zero Gromov–Witten invariant, is (rationally) homotopy equivalent to $\CP^n$, with identical first Chern class and quantum cohomology. Together with generalizations of this result to pseudo-holomorphic lines with tangency constraints, they provide a new mechanism for establishing the topological uniqueness of Liouville fillings of certain contact manifolds. We obtain filling rigidity for the ideal boundary of a class of affine varieties that admit $\CP^n$ as a compactification; notable examples include the ideal boundary of the complement of any hyperplane arrangement and $ST^*\RP^n$ for $n\ge 3$. The proofs are based on Rabinowitz–Floer homology and its variants. 
\end{abstract}
\section{Introduction}\label{s1}
Recognizing and characterizing the simplest objects in a given geometry is one of the most natural questions. 
In complex, K\"ahler, projective, or symplectic geometry, $\CP^n$ equipped with the standard complex, K\"ahler, or symplectic structure is widely regarded as the simplest closed example. 
In projective geometry, there are famous characterizations of $\CP^n$ using the Fano index by Kobayashi and Ochiai~\cite{zbMATH03411441} or the pseudo-index\footnote{The pseudo-index is the minimum of $-K_X\cdot C$ over all rational curves $C$, where $K_X$ is the canonical line bundle.} by Cho, Miyaoka, and Shepherd-Barron~\cite{zbMATH01943015} as well as Kebekus~\cite{zbMATH01803777}. 
In the symplectic category, far less is known; the only result is McDuff's theorem~\cite{zbMATH04192555} in dimension $4$: a minimal symplectic $4$-manifold is $\CP^2$ if it contains an embedded symplectic sphere of self-intersection $1$. This recognition theorem, in combination with Taubes's seminal works \cite{Taubes00}, allows Taubes to obtain the classification of symplectic structures on $\CP^2$ \cite{Taubes95}.

Both the result of Cho--Miyaoka--Shepherd-Barron and that of McDuff characterize $\CP^n$ using algebraic rational curves or symplectic spheres, which in the proofs are understood as pseudo-holomorphic spheres. 
In projective geometry, from the perspective of curves, $\CP^n$ can be characterized by the existence of lines through two generic points. Here, by a \emph{line} we mean a rational curve $C$ such that $\cL \cdot C=1$ for an ample line bundle $\cL$. 
The main result of this paper provides two generalizations of this characterization to the symplectic category and derives from them rigidity results for symplectic fillings.

\subsection{Recognition criteria}
Let $X$ be a symplectic manifold and $A\in H_2(X;\Q)$. 
Then $\GW_{X,A}(\pt,\pt)$ denotes the genus-zero Gromov--Witten invariant counting curves in class $A$ passing through two generic points of $X$.

\begin{thmx}\label{thm:main}
    Let $(X,\omega)$ be a $2n$-dimensional closed symplectic manifold with an integral symplectic form $\omega$. 
    Assume there exists a class $A\in H_2(X;\Q)$ with $$\la A,\omega \ra = 1, \quad \GW_{X,A}(\pt,\pt)\ne 0.$$ 
    Then the following hold.
    \begin{enumerate}[label=(A\arabic*)]
        \item $X$ is simply connected and rationally homotopy equivalent to $\CP^n$ and $c_1(X)=(n+1)[\omega]\in H^2(X;\Q)$.
        \item $\QH^*(X)=\QH^*(\CP^n)$ after a change of Novikov coordinate, where $\QH^*(X)$ is the small quantum cohomology.
        \item If, moreover, $\GW_{X,A}(\pt,\pt)=\pm 1$, then $X$ is homotopy equivalent to $\CP^n$. If $\GW_{X,A}(\pt,\pt)=1$, $X$ has the same genus $0$ Gromov-Witten potential as $\CP^n$.
    \end{enumerate}
\end{thmx}
When $n=2$, one can combine the argument with~\cite{zbMATH04192555} to conclude that $X$ is symplectomorphic to $\CP^2$; see \Cref{cor:McDuff}\footnote{This does not follow directly from \cite{Taubes95}, as we do not know if $\CP^2$ has a unique smooth structure.}.

In $(\CP^n,\omega_{\FS})$, rigid lines persist if one relaxes one point-constraint to a linear subspace and enforces additional tangency conditions to a local holomorphic hypersurface near the other point-constraint. This is captured by the following Gromov--Witten invariant. Let $A\in H_2(X;\Q)$ and let $B$ be a $2k$-dimensional submanifold of $X$. 
Then $\GW_{X,A}(B,\pt,k)$ is the Gromov--Witten invariant of rational curves in class $A$ with two marked points, one passing through $B$ and the other passing through a point $p$ with tangency order $k$ with respect to a local holomorphic hypersurface through $p$. We have the following generalization, which is the main force for our filling rigidity results.

\begin{thmx}\label{thm:tangency}
    Let $(X,\omega)$ be a $2n$-dimensional closed symplectic manifold with an integral symplectic form. 
    Assume there exist a homology class $A\in H_2(X;\Q)$ and a $2k$-dimensional submanifold $B\subset X$ such that the following hold:
    $$\la A,\omega \ra=1, \quad \la [B],\omega^k \ra = m \ne 0, \quad \GW_{X,A}(B,\pt,k)\ne 0.$$
    Then $X$ is simply connected, rationally homotopy equivalent to $\CP^n$, and has the same small quantum cohomology. 
    The rational homotopy equivalence can be upgraded to a homotopy equivalence if $\GW_{X,A}(B,\pt,k) = \pm m\cdot k!$.
\end{thmx}

The proofs pass to the Boothby--Wang contact manifold associated to $(X,\omega)$. 
With the help of the Rabinowitz--Floer homology developed by Cieliebak, Frauenfelder, and Oancea~\cite{zbMATH05541980,CieliebakFrauenfelderOancea2010}, adapted to a version with DGA augmentations, we then control the cohomology of the contact manifolds using the Gromov-Witten conditions. The Gromov-Witten invariant in Theorem \ref{thm:main} lifts to punctured holomorphic curves in the symplectization, which represents the death of the unit. The Gromov-Witten invariants in Theorem \ref{thm:tangency} should be interpreted in the $S^1$-equivariant Rabinowitz--Floer homology. To the best of the author's knowledge, the study of a symplectic manifold via (1)~its prequantization circle bundle (not the other way around) and (2)~such application of Rabinowitz--Floer homology is new. They yield a new mechanism for proving filling rigidities.
Uniruled conditions phrased using Gromov-Witten invariants were previously exploited by Pieloch \cite{AP25} to obtain information on $\pi_1$ and $H_2$. His methods focus on mapping degrees from moduli spaces as well as pulling cycles using universal families, which will also be used in \S \ref{s:filling} for filling rigidities, more precisely for dropping the $H_1$-injection assumption. 
The lifting perspective for computing Rabinowitz--Floer homology has appeared in the work of Bae-Kang-Kim \cite{BKK}. 

Combining Theorems~\ref{thm:main} and~\ref{thm:tangency} with the Kobayashi-Ochiai Theorem, we obtain the following.
\begin{corollaryx}\label{cor}
    If in Theorems~\ref{thm:main} or~\ref{thm:tangency}, the form $\omega$ is assumed to be K\"ahler, then $X$ is biholomorphic to $\CP^n$. 
    Moreover, if $X$ as an almost complex manifold can be equipped with a K\"ahler class $\eta$ (potentially different from $\omega$) and $n$ is odd, then $(X,\eta)$ is biholomorphic to $\CP^n$.
\end{corollaryx}

Starting from the homotopy equivalence in \Cref{thm:main}, one may try to promote it to a homeomorphism, a diffeomorphism, and finally a symplectomorphism. 
However, even in the topological category, there are infinitely many distinct homeomorphism types of manifolds that are homotopy equivalent to $\CP^n$ for $n\ge 3$, by the works of Montgomery-Yang~\cite{zbMATH03224307}, Hsiang~\cite{zbMATH03222632}, and systematically for PL structures by Sullivan~\cite{sullivan1996triangulating}. 
These structures are classified by the surgery long exact sequences, essentially by certain characteristic classes. It is interesting to investigate whether this information can be determined, and whether one can conclude that $X$ is homeomorphic to $\CP^n$.

\subsection{Applications to symplectic fillings}
The rigidity results of Theorems~\ref{thm:main} and~\ref{thm:tangency} can be turned into rigidity results for symplectic fillings via cut-and-paste constructions. 
For example, \Cref{thm:main} would provide yet another proof of the celebrated Eliashberg--Floer--McDuff theorem. 
Nowadays we have several proofs of the Eliashberg-Floer-McDuff theorem, e.g.\ by dominating the filling with holomorphic spheres in a partial compactification~\cite{zbMATH04186536} or a full compactification~\cite{zbMATH07528621}, by dominating the filling with holomorphic disks having boundary on bLobs~\cite{zbMATH06670705}, or by controlling the topology via the vanishing of symplectic (co)homology~\cite{zbMATH07367119}. 
Each perspective yields its own generalizations (see e.g.\ \cite{zbMATH06008846,zbMATH07131053,zbMATH07601288,zbMATH07757260,arXiv:2506.07525,zbMATH08152684}, \cite{zbMATH07528621}, \cite{zbMATH06670705}, and \cite{zbMATH07367119,zbMATH07673358,zbMATH07824483} respectively). 
At the level of the Eliashberg-Floer-McDuff theorem, the underlying curve is not essentially different from other proofs. 
However, the perspectives provided by Theorems~\ref{thm:main} and~\ref{thm:tangency} open a very different direction for generalization, compared to all other existing methods for filling rigidity. 
In particular, we have the following filling rigidity result for contact boundaries of a large class of affine varieties admitting $\CP^n$ as a compactification. Most of them are boundaries of non-flexible Weinstein domains, whose filling rigidity is harder to obtain. 

\begin{thmx}\label{thm:homology_unique}
    For $n\ge 3$, let $D=D_1\cup \ldots \cup D_r\subset \CP^n$ be a reduced algebraic hypersurface, not necessarily smooth. Suppose at least one irreducible component has degree less than $n$. We write $A_D=\CP^n\backslash D$ and $Y_D$ is the ideal contact boundary of $A_D$. Let $W$ be a Liouville filling of $Y_D$, such that $H_1(Y_D;\Z)\to H_1(W;\Z)$ is injective. Then:
    \begin{enumerate}[label=(D\arabic*)]
        \item\label{D1} $H_*(W;\Z)=H_*(A_D;\Z)$ as abelian groups.
        \item\label{D2} The maps $H_q(Y_D;\Z)\to H_q(W;\Z)$ are surjective and $H^q(W;\Z)\to H^q(Y_D;\Z)$ are injective. Both are isomorphisms when $q\le n-2$.
        \item\label{D3} $\pi_1(Y_D)$ normally generates $\pi_1(W)$.
        \item\label{D4} Let $d=\gcd(\deg D_1,\ldots,\deg D_r)$, the complex bundle $[TW]\oplus [\underline{\C}] = (n+1)[L_h]$ in $K^0(W)$, where $L_h$ is the complex line bundle associated to a class $h\in H^2(W;\Z)$ of torsion order $d$. In particular, all Chern classes are rationally trivial. 
    \end{enumerate}
\end{thmx}
The injective condition on the first homology is to obtain the integrality of the symplectic form in the cut-and-paste construction. Such a condition is automatic if we consider Weinstein fillings of $Y_D$. 

There are analogous results when $n=2$, however, they will not exceed the rich literature on this subject in dimension $4$.  Indeed, many filling rigidities in dimension $4$, e.g.\ \cite{zbMATH04192555,Lisca08,cap}, were obtained by recognition theorems, e.g.\ \cite{zbMATH04192555,Li06,OO96}, if a nice symplectic cap, e.g.\ a neighborhood of the ample divisor, exists. In this paper, we shall focus on the $n\ge 3$ case and use the recognition results in Theorems \ref{thm:main} and \ref{thm:tangency}.

In some cases, we can derive the condition on $H_1$ and get stronger filling rigidity results. We say an irreducible algebraic hypersurface $D\subset \C^n$ is subcritical if there exists a polynomial automorphism $\Phi:\C^n\to \C^n$, such that the defining function for $D$ is $f\circ \Phi =0$ with $f$ involving at most $n-1$ coordinates. 
\begin{thmx}\label{thm:subcritical}
    Let $D_1,\ldots, D_r$ be a sequence of subcritical hypersurfaces and $D=\cup D_i$. Then for every Liouville filling $W$ of $Y_D$, we have 
    \begin{enumerate}[label=(E\arabic*)]
        \item\label{E1} $H_1(Y_D;\Z)\to H_1(W;\Z)$ is an isomorphism and conclusions in \Cref{thm:homology_unique} hold for $W$. In particular,  $TW$ is stably trivial, hence all Chern classes are trivial. 
        \item\label{E2} If $D_1,\ldots,D_r$ are hyperplanes in $\C^n$, i.e.\ a hyperplane arrangement, then $H^*(W;\Z)\hookrightarrow H^*(Y_D;\Z)$ is isomorphic to $H^*(A_D;\Z)\hookrightarrow H^*(Y_D;\Z)$ as maps of rings. In particular, $H^*(W;\Z)$ is isomorphic to the Orlik–Solomon algebra \cite[Definition 3.45]{zbMATH00052944} of the hyperplane arrangement. 
    \end{enumerate}
\end{thmx}
A special case of hyperplane arrangements is when $D_1,\ldots,D_n$ are the coordinate hyperplanes, then $Y_D=ST^*T^n$, whose Liouville filling rigidity up to diffeomorphism was obtained by Bowden-Gironella-Moreno \cite{zbMATH07601288} and Geiges-Kwon-Zehmisch \cite{zbMATH07757260}. In particular, $A_D$ is generally a critical Weinstein domain.

Another important family of filling rigidity comes from considering the quadratic hypersurface in $\CP^n$ for $n\ge 3$.  In this case, $Y_D$ is $ST^*\RP^n$ equipped with the canonical contact structure, and we can also drop the condition on $H_1$ with a careful analysis of holomorphic spheres and their universal families to probe the first homology.

\begin{thmx}\label{thm:RP}
    For $n\ge 3$, let $W$ be a Liouville filling of $Y=ST^*\RP^n$ equipped with the canonical contact structure. Then we have 
    \begin{enumerate}[label=(F\arabic*)]
        \item\label{F1} $W$ has the same integral homology group and integral cohomology ring as $\RP^n$. $\pi_1(Y)\to \pi_1(W)$ is an isomorphism. $[TW]\oplus [\underline{\C}]=(n+1)[L_h]\in K^0(W)$ where $L_h$ is the complex line bundle associated to the non-trivial class $h\in H^2(W;\Z)$.
        \item\label{F2} If $n$ is odd, then $W$ is diffeomorphic to $DT^*\RP^n$.
        \item\label{F3} If $n$ is even and $W$ is Weinstein, then $W$ is diffeomorphic to $DT^*\RP^n$. 
    \end{enumerate}
\end{thmx}
Note that \Cref{thm:RP} is optimal. Indeed, rigidity for fillings of $ST^*\RP^n$ fails when $n=2$, since $ST^*\RP^2=(L(4,1),\xi_{\std})$ admits two different Liouville fillings with distinct rational cohomology groups~\cite{zbMATH04192555}. This justifies the necessity of the $k<n$ condition in \Cref{thm:homology_unique}, which is used in an essential way to invoke \Cref{thm:tangency}.

\subsection*{Organization of the paper}
In \S\ref{s:RFH} we review the Rabinowitz--Floer homology developed by Cieliebak, Frauenfelder, and Oancea, and introduce the version used in this paper. 
We prove \Cref{thm:main} in \S\ref{s:top} and \Cref{thm:tangency} in \S\ref{s4}.  
Finally, in \S\ref{s:filling} we discuss applications to symplectic fillings and prove Theorems~\ref{thm:homology_unique}-\ref{thm:RP}, as well as other refinements.

\subsection*{Acknowledgments}
The author would like to thank Wenhao Ou for explaining the algebro-geometric proof of Corollary~\ref{cor} under the assumptions of \Cref{thm:main}, Xiaobo Liu for explaining the uniqueness of the genus $0$ Gromov-Witten potential,  Kai Cieliebak, and Jungsoo Kang for helpful conversations. 
The author is supported by the National Key R\&D Program of China under Grant No.\ 2023YFA1010500, the National Natural Science Foundation of China under Grant No.\ 12625108, 12288201 and 12231010.

\subsection*{AI statement}
I acknowledge the use of GPT-6 Astra for proofreading and literature search. In particular, a gap was found by AI in the proof of dropping the $H_1$-injective condition for divisors of prime degree at most $(n+1)/2$, which was announced in some talks.

\section{Rabinowitz Floer homology}\label{s:RFH}
In this part, we recall the basics of Rabinowitz Floer homology developed by Cieliebak, Frauenfelder, and Oancea \cite{zbMATH05541980,CieliebakFrauenfelderOancea2010}. A good summary would be \cite{CieliebakOancea2018}. As the literature on Rabinowitz Floer homology, unlike symplectic (co)homology, is rather consistent with the convention of homology or cohomology, we will work with the same convention as in \cite{CieliebakOancea2018}. In particular, symplectic homology in this paper will be the symplectic cohomology (after a grading modification) in almost all of the author's other papers \cite{zbMATH07367119,zbMATH07673358,zbMATH07824483} as well as those in \cite{zbMATH05504309}. All coefficients are taken in a field $\mathbb K$; this ensures exactness of both the inverse limit over finite‑dimensional vector spaces and tensor products.

\subsection{Symplectic homology of filled cobordism}
\label{subsec:setup}
Let  $(W^{2n},\lambda)$ be an exact cobordism, namely
\begin{enumerate}
    \item $\omega=\rd\lambda$ is a symplectic form;
    \item The Liouville vector field $Z$ defined by $\iota_Z \rd \lambda=\lambda$ is transverse to the boundary of $W$;
    \item The boundary $\partial W$ splits into the convex/positive part $\partial_+W$ and the concave/negative part $\partial_-W$, depending on whether $Z$ is outward or inward pointing. 
    \item We use $\alpha_{\pm}$ to denote the restrictions of $\lambda$ to $\partial_{\pm}W$, which are contact forms on the boundary.
\end{enumerate}
If $\partial_-W=\emptyset$, $(W,\lambda)$ is called a Liouville filling of $(\partial_+W,\alpha_+)$. Throughout the paper, we will consider only connected symplectic cobordisms/manifolds.
 
\begin{definition}
    A filled cobordism $(W,F,\lambda)$ is a Liouville cobordism $(W,\lambda)$, along with a Liouville filling $(F,\lambda_F)$ of $(\partial_-W,\alpha_-)$. We can glue $F$ to $W$, to obtain a Liouville filling $W_F:=F\circ W$ of $(\partial_+W,\alpha_+)$. 
\end{definition}
We use $(\widehat{W},\widehat{\lambda})$ to denote the completion, namely the gluing of $W$ with two cylindrical ends
$$\left((0,1)_r\times \partial_-W \cup W \cup (1,\infty)_r\times \partial_+W, \quad r\alpha_-\cup \lambda \cup r\alpha_+\right).$$
Then $\widehat{W}_F$ contains $\widehat{W}$. 

On an exact symplectic manifold $(W,\lambda)$, given a Hamiltonian $H$, the convention for Hamiltonian vector field $X_H$ is $\iota_{X_H}\rd \lambda=-\rd H$ and the symplectic action of a loop $x:S^1\to W$ for an $S^1$-dependent Hamiltonian $H_t$ is 
\[
\cA_H(x)=\int_{S^1}x^*\lambda-\int_{S^1}H_t(x(t))\rd t
\]
On a contact manifold $(M,\alpha)$, the Reeb vector field $R_\alpha$ is defined by $i_{R_\alpha}d\alpha=0$, $\alpha(R_\alpha)=1$. Its closed orbits are called Reeb orbits; the set of their periods is $\Spec(M,\alpha)$, a closed subset of $(0,\infty)$. For a Hamiltonian $H(r,x)=h(r)$ on the symplectization $((0,\infty)\times M, r\alpha)$, the Hamiltonian vector field is $X_H = h'(r)R_\alpha$, so $1$-periodic orbits on level $\{r\}\times M$ correspond to Reeb orbits of period $h'(r)$. We use $\Spec(W,\lambda)$ to denote $\Spec(\partial_+W, \alpha_+)\cup \{0\} \cup (-\Spec(\partial_-W,\alpha_-))$.

We use $\cH(W;F)$ to denote the set of admissible Hamiltonians on $\widehat{W}_F$, i.e.\ $S^1$-dependent Hamiltonians with the following properties:
\begin{enumerate}
    \item $H_t<0$ on $W$.
    \item $H_t=ar+c$ for $r>1+\delta$ on $(1,\infty)_r\times \partial_+W \subset \widehat{W}_F$, for $a\notin \Spec(\partial_+W,\alpha_+)$.
    \item $H_t$ has non-degenerate Hamiltonian orbits
\end{enumerate}
$\cH(W;F)$ is a partial order set, with $H_t\le G_t$ if $H_t(x)\le G_t(x)$ for all $t\in S^1$ and $x\in \widehat{W}_F$.

Given an admissible Hamiltonian $H_t$ and a compatible $J_t$ that is cylindrical (i.e.\ $J_t$ preserves the contact structure, $J(r\partial_r)=R$ and is independent of $r$) at the end of $\widehat{W}_F$, we consider the Gromov-Floer compactification $\cM_{H_t,J_t}(x_+,x_-)$ of the moduli space
$$\widetilde{\cM}_{H_t,J_t}(x_+,x_-)=\left\{u:\R_s\times S^1_t\to \widehat{W}_F\left| \partial_su+J_t(\partial_tu-X_{H_t})=0, \lim_{s\to \pm \infty }u = x_{\pm} \right.\right\}/\R,$$
where $x_{\pm}\in \cP(H)$ are two Hamiltonian orbits of $H_t$. In the following, we will suppress the subscript $t$ in $H_t$ and $J_t$, and suppress $J_t$ in $\cM_{H_t,J_t}$. 
Then we have: 
\begin{enumerate}
    \item For $u\in \widetilde{\cM}_{H}(x_+,x_-)$, we have $\cA_{H}(x_+)-\cA_{H}(x_-)=\int_{\R\times S^1} |\partial_su|^2\rd s\rd t\ge 0$, computed using the metric induced from $J$. 
    \item For generic choice of $J$, $\widetilde{\cM}_H(x_+,x_-)$ is a manifold of dimension $\mu_{\CZ}(x_+,u)-\mu_{\CZ}(x_-,u)-1$, where the Conley-Zehnder index is understood with respect to a symplectic trivialization of $u^*T\widehat{W}_F$. $\mu_{\CZ}(x_+,u)-\mu_{\CZ}(x_-,u)-1$ is referred to as the virtual dimension and only depends on the homotopy class of $u$. We use $\widetilde{\cM}^k_H(x_+,x_-)$ to denote the component with virtual dimension $k$. Without the reference of $u$, the Conley-Zehnder index $\mu_{\CZ}(x)$ is only well-defined in $\Z/2$ unless we fix a trivialization of $\det_{\C}T\widehat{W}_F$.
    \item If we assume the moduli spaces of expected dimension at most $1$ are cut out transversely, then $\cM^0_H(x_+,x_-)=\widetilde{\cM}^0_H(x_+,x_-)$ and $\partial \cM^1_H(x_+,x_-)=\cup_y\cM^{0}_H(x_+,y)\times \cM^0_H(y,x_-)$.
\end{enumerate}

Let $\K$ be a field, for $a<b \notin \Spec(W,\lambda)$, we define $$\FC_*^{<b}(H)=\displaystyle \bigoplus_{\substack{x\in \cP(H)\\ \cA_H(x)<b}}\K \cdot x, \quad \FC_*^{(a,b)}(H)=\FC^{<b}_*(H)/\FC^{<a}_*(H)\footnote{Strictly speaking, this should be denoted by $\FC_*^{[a,b)}(H)$. However, after taking the direct limit for $H$, with the assumption that $a,b\notin \Spec(W,\lambda)$, the chain complex/homology does only see generators with action window $(a,b)$. We use such notation to be consistent with \cite{CieliebakOancea2018}.}$$
graded by the Conley-Zehnder indices in $\Z/2$ (in $\Z$, if we fix a trivialization of $\det_{\C}T\widehat{W}_F$), with a differential $\partial:\FC_*^{(a,b)}(H)\to \FC_{*-1}^{(a,b)}(H)$ defined by 
$$\partial(x_+)=\sum \#\cM^0_H(x_+,x_-)\cdot x_-.$$
Then we define $\FH^{(a,b)}_*(H)=H_*(\FC^{(a,b)}_*(H),\partial)$. Now if $H_1\le H_2$, by counting solutions to $\partial_su+J(\partial_tu-X_{H_s})=0$ with $\rd H_s/\rd s\le 0$ , $H_s=H_2$ for $s\ll 0$ and $H_s=H_1$ for $s\gg 0$, we get a continuous map $\FH^{(a,b)}_*(H_1)\to \FH^{(a,b)}_*(H_2)$. The \emph{filtered symplectic homology} of $W$ (with respect to filling $F$) is
\[
\SH_*^{(a,b)}(W;F) := \varinjlim_{H\in\cH(W;F)} \FH_*^{(a,b)}(H).
\]

To get a clear picture of generators for the filtered symplectic homology. We consider the cofinal subset $\cH_V(W;F)$ of V-shaped Hamiltonians that consists of perturbations of the following 
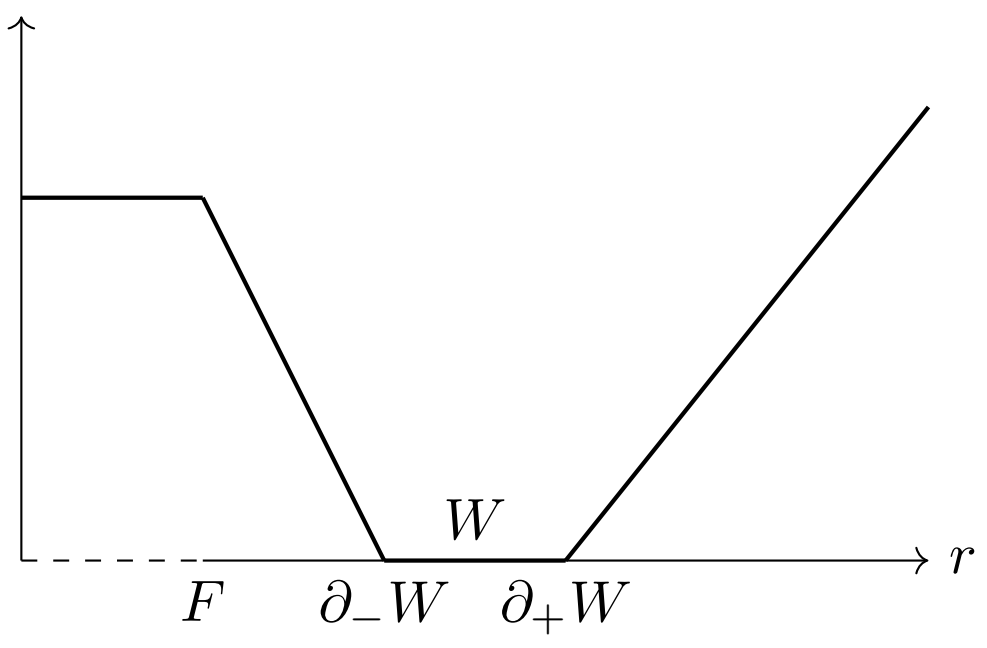
\begin{figure}[H]
\centering
\begin{tikzpicture}[scale=1.2]
\draw[->] (1,0) -- (5,0) node[right] {$r$};
\draw[dashed] (0,0) -- (1,0) node[below] {$F$};
\draw[->] (0,0) -- (0,3);
\draw[thick] (0,2) -- (1,2);
\draw[thick] (1,2) -- (2,0);
\draw[thick] (2,0) -- (3,0) node[midway,above] {$W$};
\draw[thick] (3,0) -- (5,2.5);
\draw  (2,0) node[below] {$\bd_-W$};
\draw (3,0) node[below] {$\bd_+W$};
\end{tikzpicture}
\caption{V‑shaped Hamiltonian, cf. \cite[Figure~8]{CieliebakOancea2018}.}
\label{fig:Vshape}
\end{figure}
\noindent Namely, $H$ is zero on $W$ and is linear slope $\tau>0$ on $(1,\infty)\times \partial_+W$ and is of slope $-\mu<0$ on $(1-\delta,1)_r\times \partial_-W \subset \widehat{W}_F$ and is $\mu\delta$ on the rest of domain contained in $F$. The perturbation is $C^2$-small on the domain where $H=0$, producing only constant Hamiltonian orbits, and will be in the form of functions that only depend on $r$ on the part where $H$ is linear. In particular, Hamiltonian orbits in $S^1$-families correspond to Reeb orbits. One may perturb further to break the $S^1$-symmetry. There are the following five types of Hamiltonian orbits:
\begin{enumerate}
    \item Type $F$: constant orbits on the subset of $F$ where $H=\mu\delta$ before perturbation, their symplectic action is $\approx -\mu\delta$. 
    \item Type $F_-$: non-constant orbits near $\{1-\delta\}\times \partial_-W$; They, in pairs, correspond to \emph{negatively} parametrized Reeb orbits of $(\partial_-W,\alpha_-)$ of period  at most $\mu$, whose symplectic action is approximately $-\mu\delta-(1-\delta)T$. Here $T$ is the period of the corresponding Reeb orbit $\gamma$. In particular, their action is smaller than $-\mu\delta$. The Conley-Zehnder indices of the two orbits are $-\mu_{\CZ}(\gamma),-\mu_{\CZ}(\gamma)-1$, where $\mu_{\CZ}(\gamma)$ is the Conley-Zehnder index of the Reeb orbit $\gamma$ (defined in $\Z$ after we fix a trivalization of $\gamma^*\det_{\C}\xi$). Here, the Conley-Zehnder indices of the Hamiltonian orbits are computed using a symplectic trivialization of the contact structure $\gamma^*\xi$, with the symplectic normal canonically framed by the Reeb vector field. 
    \item Type $I_-$: non-constant orbits near $\{1\}\times \partial_-W$; They, in pairs, correspond to \emph{negatively} parametrized Reeb orbits of $(\partial_-W,\alpha_-)$ of period  at most $\mu$, whose symplectic action is approximately $-T$, the negative period of the corresponding Reeb orbit $\gamma$.  The Conley-Zehnder indices of the two orbits are $-\mu_{\CZ}(\gamma)+1,-\mu_{\CZ}(\gamma)$.
    \item Type $I_0$: constant orbits on $W$, their symplectic action is $\approx 0$.
    \item Type $I_+$: non-constant orbits near $\{1\}\times \partial_+W$; They, in pairs, correspond to \emph{positively} parametrized Reeb orbits of $(\partial_+W,\alpha_+)$ of period  at most $\tau$, whose symplectic action is approximately the period $T$ of the corresponding Reeb orbit $\gamma$. The Conley-Zehnder indices of the two orbits are $\mu_{\CZ}(\gamma),\mu_{\CZ}(\gamma)+1$.
\end{enumerate}
Then the direct limit for $H\in \cH(W;F)$ is equivalent to the direct limit of $\tau,\mu\to \infty$ in $\cH_V(W;F)$. For $a,b\notin \Spec(W,\lambda)$, we can compute $SH^{(a,b)}_*(W;F)$ using a Hamiltonian $H\in \cH_V(W;F)$ with $\tau>\max\{b,0\}$ and $\mu>\max\{-a/\delta,0\}$, then orbits of type $F$ and $F_-$ are filtered out. The chain complex only have generators from orbits of type $I_+$, such that the corresponding Reeb orbits on $\partial_+W$ have period in $(a,b)$, orbits of type $I_-$, such that the corresponding Reeb orbits on $\partial_-W$ have period in $(-b,-a)$, and orbits of type $I_0$ if $0\in (a,b)$.

\begin{definition}[{\cite[Definition 2.8]{CieliebakOancea2018}}]
We define six versions of symplectic homology groups of $W$:
\[
\begin{array}{rclr}
\SH_*(W;F) &:= & \displaystyle\varinjlim_{b\to\infty}\varprojlim_{a\to-\infty}\SH_*^{(a,b)}(W;F) & \text{(full symplectic homology)},\\
\SH_*^{>0}(W;F) &:= &  \displaystyle\varinjlim_{b\to\infty}\varprojlim_{a\searrow0}\SH_*^{(a,b)}(W;F) &\text{(positive symplectic homology)},\\
\SH_*^{\ge0}(W;F) &:= & \displaystyle\varinjlim_{b\to\infty}\varinjlim_{a\nearrow0}\SH_*^{(a,b)}(W;F) &\text{(nonnegative symplectic homology)},\\
\SH_*^{=0}(W;F) &:= &  \displaystyle\varprojlim_{b\searrow0}\varinjlim_{a\nearrow0}\SH_*^{(a,b)}(W;F) &\text{(zero‑level symplectic homology)},\\
\SH_*^{\le0}(W;F) &:= &  \displaystyle\varprojlim_{b\searrow0}\varprojlim_{a\to-\infty}\SH_*^{(a,b)}(W;F) &\text{(nonpositive symplectic homology)},\\
\SH_*^{<0}(W;F) &:=&  \displaystyle \varinjlim_{b\nearrow0}\varprojlim_{a\to-\infty}\SH_*^{(a,b)}(W;F) &\text{(negative symplectic homology)}.
\end{array}
\]
\end{definition}

The Rabinowitz Floer homology is then a special case, namely the symplectic homology of the trivial cobordism of a contact manifold $Y$, along with a Liouville filling $F$, i.e.
$$\RFH_*(Y;F):=\SH_*([1,2]_r\times Y;F)$$
We use $\RFH^{\heartsuit}_*(Y;F)$ for $\heartsuit\in \{\emptyset, >0,<0,\ge 0, \le 0,=0\}$ to denote the corresponding version. 
\subsection{Properties of Symplectic homology}. 
We now list the properties of symplectic homology that will be relevant for our purpose. 
\begin{proposition}\label{prop:basic_property}
We have the following properties:
\begin{enumerate}
    \item\label{bp1}  The six versions of symplectic homology are independent of the Liouville forms on $W$, but would depend on the filling $F$ up to Liouville homotopy \cite[Proposition 2.9]{CieliebakOancea2018}. The filtered version $\SH^{(a,b)}_*(W;F)$ depends on the Liouville form of $W$.  
    \item\label{bp2} $\SH^{=0}_*(W;F)=H^{n-*}(W)$ \cite[Proposition 2.10]{CieliebakOancea2018}, where $n=\dim_{\C}W$.
    \item\label{bp3} There are various long exact sequences \cite[Proposition 2.12]{CieliebakOancea2018}, e.g.
    $$\ldots \to \SH^{\le 0}_*(W;F)\to \SH_*(W;F)\to \SH^{>0}_*(W;F)\to \SH^{\le 0}_{*-1}(W;F)\to \ldots $$
    Moreover, if we decompose the complex into three pieces instead of two, we get a spectral sequence converging to $SH_*(W;F)$, with the first page $\SH^{<0}_*(W;F)\oplus \SH^{=0}_*(W;F)\oplus \SH^{>0}_*(W;F)$. The differential on the first page is the connecting map 
    \begin{equation}\label{eqn:delta}
        \delta:\SH^{>0}_*(W;F)\to \SH^{=0}_{*-1}(W;F), \quad \delta^{\vee}: \SH^{=0}_*(W;F)\to \SH^{<0}_{*-1}(W;F).
    \end{equation}
    from the tautological long exact sequences relating $\heartsuit\in \{\ge 0,>0,=0\}$ and $\heartsuit\in \{\le 0, <0,=0\}$
\end{enumerate}
\end{proposition}
\begin{proof}
    The three-piece claim is not present in \cite{CieliebakOancea2018}, but it is identical to the proof of \cite[Proposition 2.12]{CieliebakOancea2018}. Namely for $a<0$ and $b>0$, we have a spectral sequence convergent to $\SH^{(a,b)}_*(W;F)$ with first page $\SH^{(a,-\epsilon)}_*(W;F)\oplus \SH^{(-\epsilon,\epsilon)}_*(W;F)\oplus \SH^{(\epsilon,b)}_*(W;F)$ for $0<\epsilon<\min(\Spec(\partial_+W,\alpha_+)\cup \Spec(\partial_-W,\alpha_-))$ with differential the filtered version of $\delta$ and $\delta^{\vee}$. The spectral sequence is compatible with increasing $b$ and decreasing $a$, since the inverse limit for $a\to -\infty$ is exact as $SH^{(a,b)}_*(W;F)$ is a finite dimensional $\K$-space. Applying inverse limit for $a\to -\infty$ and then direct limit $b\to \infty$, we get the claimed spectral sequence computing $\SH_*(W;F)$.
\end{proof}

\begin{proposition}[{\cite[Theorem 9.4]{CieliebakOancea2018}}]\label{prop:dual}
We have a duality $\RFH^{<0}_*(Y;F)=(\RFH^{>0}_*(Y;F))^{\vee}[-1]$, and $\delta^{\vee}$ from \eqref{eqn:delta} is the dual of $\delta$, using the Poincar\'e duality on $Y$, i.e.\ $(\RFH^{=0}_*(Y))^{\vee}=(H^{n-*}(Y))^{\vee}\stackrel{\PD}{=}H^{n-1+*}(Y)=\RFH^{=0}_{-*+1}(Y)$. 
\end{proposition}
\begin{proof}
    This follows from \cite[Theorem 9.4]{CieliebakOancea2018} for $\heartsuit\in \{\le 0,<0\}$ and the fact that $\SH^*_{\heartsuit}(Y;F)=(\SH_*^{-\heartsuit}(Y;F))^{\vee}$ (the dual of direct limit is the inverse limit of the dual). 
\end{proof}
\begin{remark}
    The version of \Cref{prop:dual} for $\heartsuit\in\{\ge 0, >0\}$ does not hold in the category of vector spaces, but holds in the category of Tate vector spaces \cite{arXiv:2407.21741}. 
\end{remark}
The symplectic homology of cobordism also has a TQFT structure \cite[\S 10.1]{CieliebakOancea2018}, however, in order to maintain the maximum principle, there are constraints on the spectrum window. 
\begin{proposition}[{\cite[(64)]{CieliebakOancea2018}}]\label{prop:product}
    We have a pair of pants product
    $$SH^{(a,b)}_*(W;F)\otimes SH^{(a,b)}_*(W;F)\to SH^{(a+b,2b)}_*(W;F)[-n]$$
    compatible with the maps in the double directed system of $a,b$. As corollaries, we have $\SH^{\le 0}_*(W;F)$, $\SH^{<0}_*(W;F)$, $\SH_*(W;F)$,$\SH^{=0}_*(W;F)$ are super-commutative algebras, where $\SH^{\le 0}_*(W;F)$, $\SH_*(W;F)$, $\SH^{=0}_*(W;F)$ are unital. The product on $\SH^{=0}_*(W;F)$ is the cup product on $H^*(W)$.
\end{proposition}
By the discussion before \cite[Proposition 10.1]{CieliebakOancea2018}, $\SH^{\le 0}_*(W;F)$,$\SH_*(W;F)$,$\SH^{=0}_*(W;F)$ are unital, but $\SH^{<0}_*(W;F)$ is not. Moreover, by \cite[(63)]{CieliebakOancea2018}, we can deduce the product is defined on
$$\SH^{<b}_*(W;F)\otimes \SH^{<b'}_*(W;F)\to \SH^{<b+b'}_*(W;F).$$

\begin{proposition}\label{prop:kill}
    Let $x\in \SH^{\le 0}_*(W;F)$ with $x_0$ the image of $x$ in $\SH^{=0}_*(W;F)$ from tautological exact sequence. Assume $x_0$ is in the form of $1+A$ for $A\in \oplus_{k>0} H^k(W)$ under the identification $\SH^{=0}_*(W;F)=H^{n-*}(W)$, then $x$ is invertible in $\SH^{\le 0}_*(W;F)$.   
\end{proposition}
\begin{proof}
    Let $e$ be the identity in $\SH^{\le 0}_*(W;F)$, then the image $e_0$ of $e$ in $\SH^{=0}(W;F)$ is $1\in H^0(W)$ under the identification $\SH^{=0}_*(W;F)=H^{n-*}(W)$. Let $\eta=e-x$, then $\eta_0=-A$ by assumption. By \Cref{prop:product}, we $\eta_0^{2n}=0$, as a consequence, $\eta^{2n}=\iota(\eta_+)$ for $\eta_+\in SH^{<b}_*(W;F)$ for a fixed $b<0$, where $\iota:\SH^{<b}_*(W;F)\to \SH^{\le 0}_*(W;F)$. Hence $\eta^k_+$ is from $\SH^{<kb}(W;F)$. 
    Then we have $\sum_{i=0}^{\infty}\eta^i$ is well defined in $\SH^{\le 0}_*(W;F)$. Indeed, it can be written as 
    $$\sum_{j=0}^{N-1}\eta^j+\sum_{i\ge 1}^{\infty}\sum_{j=0}^{2n-1}\iota(\eta^{j}\eta_+^{i})$$
    which a finite sum in each $SH^{(a,0]}_*(W;F)$ for any $a<0$. Since $x^{-1}=(e-\eta)^{-1}=\sum_{i=0}^\infty \eta^i$, we have $x$ is invertible in  $\SH^{\le 0}_*(W;F)$.   
\end{proof}

\begin{corollary}\label{cor:surj}
    Let $x\in \SH^{(0,b)}_*(W;F)$, assume that $\pi_0\circ \delta(x)\ne 0$ for $\pi_0:\SH^{=0}_*(W;F)\simeq H^{n-*}(W)\to H^0(W)$ and $\delta:\SH_*^{(0,b)}(W;F)\to \SH^{=0}_{*-1}(W;F)$, then $SH^{(0,b)}_*(W;F)\to SH^{\le 0}_{*-1}(W;F)$ is surjective.
\end{corollary}
\begin{proof}
Let $u$ be the image of $x$ under the map $\SH_*^{(0,b)}(W;F) \to \SH^{\le 0}_{*-1}(W;F)$. By assumption and \Cref{prop:kill}, $u$ is a unit in $\SH^{\le 0}_*(W;F)$. Then from the tautological exact sequence, we have the unit $u$ in $\SH^{\le 0}_*(W;F)$ is mapped to zero in $\SH^{<b}_*(W;F)$. Then any $\SH^{\le 0}_*(W;F)\to \SH^{<b}_*(W;F)$ is zero as any element $x$ in $\SH^{\le 0}_*(W;F)$ can be written as $u\cdot u^{-1}x$, and $\SH^{\le 0}_*(W;F)\to \SH^{<b}_*(W;F)$ is compatible with the product with $\SH^{\le 0}_*(W;F)$. Then, from the long exact sequence, we have $SH^{(0,b)}_*(W;F)\to SH^{\le 0}_{*-1}(W;F)$ is surjective. 
\end{proof}

\subsection{Symplectic homology of augmented cobordism}
Ideally, we want to define symplectic homology just for the completed cobordism $\widehat{W}$. There are two directions of non-compactness, one from the positive end of $\widehat{W}$, for which we have the maximum principle, the other from the negative end $\widehat{W}$, for which we need to work with the filled version $\widehat{W}_F$ to obtain a tautological $ C^0$-bound. However, one can work with the SFT compactification by including degeneration from the negative end of $\widehat{W}$ and use the algebraic replacement of $F$, i.e.\ a DGA augmentation to the contact DGA of $\partial_-W$ to define a symplectic cohomology of cobordisms. Similar constructions (the $\SH^{>0}$ version) were already used in \cite{AOT}, \cite[\S 9.5]{CieliebakOancea2018} and \cite[\S 5.2]{BGMZ} (both for trivial augmentations). 

\subsubsection{Contact DGAs}
We briefly recall contact homology \cite{Bao-Honda,Pardon:Contact} and augmentations. Let $(Y,\xi)$ be a $(2n-1)$-dimensional contact manifold. We fix a contact form $\alpha$ such that all Reeb orbits are non-degenerate. For each good Reeb orbit $\gamma$, we associate it with a formal variable $q_{\gamma}$. We use $V$ to denote the $\Z/2$ graded $\K$ vector space generated by $q_{\gamma}$, with grading $|q_{\gamma}|=\mu_{\CZ}(\gamma)+n-3$ and the characteristic of the field $\K$ is zero. By counting rigid punctured holomorphic curves with one positive puncture and multiple negative punctures in the symplectization $\widehat{Y}=\R\times Y$, we can get degree $-1$ operations
$$\partial_k:V\to S^kV=(\otimes^k V)/\Sigma_k$$
such that $\partial = \sum_{k=0}^\infty \partial_k$ on $SV=\oplus_{k=0}^\infty S^kV$ defined by the Leibniz rule 
$$\partial(ab)=\partial(a)b+(-1)^{|a|} a\partial(b)$$
determines a differential graded algebra (DGA) $(SV,\partial)$. 
In a slightly more precise manner, we fix an $\rd(e^s\alpha)$-compatible $s$-invariant almost complex structure on the completion $\widehat{Y}=\R_s\times Y$ such that $\alpha \circ J = \rd s$, i.e., cylindrically convex. For each Reeb orbit, we fix a base point on the image. Let $\gamma$ be a good Reeb orbit and $\Gamma=\{\gamma_1,\ldots,\gamma_m\} $ be an ordered multiset (i.e., a set allowing duplicates) of good Reeb orbits; we define $\cM_{\gamma,\Gamma}$ to be the compactified moduli space
\begin{equation}
    \overline{
\left\{ 
\begin{array}{c}
     u: \CP^1\backslash \{z_+,z_{1,-},\ldots,z_{m,-}\} \to \R_s\times Y\\ 
     z_+,z_{1,-},\ldots,z_{m,-}\in \CP^1 \text{ pairwise different} \\
     \theta_{+} \in ST_{z_+}\CP^1, \theta_{i,-} \in ST_{z_{i,-}}\CP^1
\end{array}\left|
\begin{array}{c}
     (\rd u)^{0,1}=0,\\
     u(z_+) = \gamma,\\
     u(z_{i,-})=\gamma_{i}.
\end{array}
\right.
\right\}/\mathrm{Aut}\times \R    
    }
\end{equation}
Here $u(z_+) = \gamma$ is a shorthand for $u$ converging to $(\infty,\gamma)$ near the positive puncture $z_+$ on $\CP^1$ such that the asymptotic marker $\theta_+$ is mapped to the base point on the image of $\gamma$; similarly $u(z_{i,-}) = \gamma_i$ is a shorthand for $u$ converging to $(-\infty,\gamma_i)$ near the negative puncture $z_{i,-}$ with a similar condition on the asymptotic marker $\theta_{i,-}$. The automorphism group consists of biholomorphisms of $\CP^1$ intertwining all data, i.e., $u,z_+,z_{1,-},\ldots,z_{m,-}, \theta_+, \theta_{1,-},\ldots,\theta_{m,-}$. And the $\R$-action is from the ambient $s$-translation on $\R_s\times Y$. Then the operator $\partial_k$ is defined by 
$$\partial_k(q_{\gamma}) = \sum_{[\Gamma],|\Gamma|=k} \frac{1}{\mu_{\Gamma}\kappa_{\Gamma}} \# \cM_{\gamma,\Gamma}q^{\Gamma},$$
where $q^{\Gamma}:=q_{\gamma_1}\cdot \ldots \cdot q_{\gamma_m}$. The sum is over all (un-ordered) multiset $[\Gamma]$ with $k$ elements such that $\cM(\gamma,\Gamma)$ has virtual dimension $0$. $\Gamma$ is an ordered representative, $\kappa_{\Gamma}$ is the product of multiplicities of orbits in $\Gamma$ and $\mu_{\Gamma}$ is the size of the isotropy group of $\Gamma$ from reordering, i.e., if $\Gamma=\{\underbrace{\gamma_1,\ldots,\gamma_1}_{i_1},\ldots, \underbrace{\gamma_k,\ldots,\gamma_k}_{i_k}\}$ for $\gamma_i\ne \gamma_j$ if $i\ne j$, then $\mu_{\Gamma}=i_1!\cdot \ldots \cdot i_k!$. Finally, the counting $\# \cM_{\gamma,\Gamma}$ requires fixing a virtual setup \cite{Bao-Honda,Pardon:Contact} and extra auxiliary choices/perturbations therein.

We use $\CC(Y,\alpha,\Theta)$ to denote this DGA as it depends on the contact form $\alpha$ as well as a set of auxiliary data $\Theta$ depending on the choice of virtual setup. The homology of $\CC(Y,\alpha,\Theta)$ is called the contact homology of $Y$, which is a contact invariant \cite{Bao-Honda,Pardon:Contact}\footnote{Strictly speaking, \cite{Bao-Honda} did not give a construction of the full contact homology algebra $\CC(Y,\alpha,\Theta)$, but only a filtered version $\CC^{L_i}(Y,\alpha_i,\Theta_i)$ by Reeb orbits of period at most $L_i$ of a sequence of contact forms $\alpha_i$ that are only $C^0$ close to each other. They form a direct system only on the homology level.}, i.e., it only depends on $\xi$. In this paper, we will use Pardon's framework \cite{Pardon:Contact}.

An augmentation $\epsilon$ is a DGA map from $(SV,\partial)$ to $(\K,0)$, i.e., a map $\epsilon:SV\to \K$ such that $\sum \epsilon^k \circ \partial_k:V\to \K$ is zero. Given an augmentation $\epsilon$, we have an algebra isomorphism $F_{-\epsilon}:SV\to SV$ determined by $F_{-\epsilon}(x)=x-\epsilon(x)$. The inverse of $F_{-\epsilon}$ is $F_{\epsilon}$. Under such change of coordinate, $F_{\epsilon}\circ \partial \circ F_{-\epsilon}$ is again induced from a sequence of operations $\partial_{\epsilon,k}:V\to S^kV$ with the significance that $\partial_{\epsilon,0}=0$. As a consequence, we have $\partial_{\epsilon,1}$ is a differential on $V$. Then we can define the linearized contact homology w.r.t. $\epsilon$ by $\LCH_{*}(Y,\epsilon)=H_*(V,\partial_{\epsilon,1})$. The linearized contact homology certainly depends on the augmentation, but the set of linearized contact homology after enumerating through all augmentations is a contact invariant by \cite{chaidez2024contact}. As the contact DGA of $\emptyset$ is the ground field $\K$, counting rigid holomorphic planes in an exact filling $W$, plus suitable virtual/perturbation data, yields an augmentation to the contact DGA by \cite{Pardon:Contact}.

\subsubsection{Constructions of Symplectic homology of augmented cobordism}\label{ss:SH_aug} 
Symplectic (co)homology of augmentations appeared in \cite{AOT,BGMZ}, which are the $\ge 0, >0$ versions of the symplectic homology of augmented cobordism in this paper. Such ideas were originally due to Bourgeois-Oancea \cite{BO09}.

Let $(W,\lambda)$ be an exact cobordism. Let $\alpha_0$ be a contact form on $Y:=\partial_-W$, such that $\alpha_0<\alpha_-$ and all Reeb orbits of $\alpha_0$ are non-degenerate. We consider the completion $\widehat{W}$, which contains $((0,1)_{r'}\times Y, r'\alpha_0)$ inside $((0,1)_{r}\times \partial_-W,r\alpha_-)$ for $r=fr'$ where $\alpha_0=f\alpha_-$. We use $\cH(W)$ to denote the admissible Hamiltonians on $\widehat{W}$, namely $H$ with following properties:
\begin{enumerate}
    \item $H<0$ on $W$.
    \item $H=ar+c$ for $r>1+\delta$ on $(1,\infty)_r\times \partial_+W\subset \widehat{W}$ and $a\notin \Spec(\partial_+W,\alpha_+)$.
    \item $H$ is constant on a closed set $K$ containing the negative end $((0,1)_{r'}\times Y, r'\alpha_0)$ of $\widehat{W}$.
    \item Hamiltonian orbits on the complement of $K$ are non-degenerate. 
\end{enumerate}
Similarly, in the poset $\cH(W)$, we have a cofinal subset $\cH_V(W)$ of V-shaped Hamiltonian. Namely, perturbations to Hamiltonians that are zero on $W$, linear with slope $\tau>0$ on $(1,\infty)\times \partial_+W$, linear with slope $-\mu<0$ on $(1-\delta,1)_r\times \partial_-W$ and is $\mu\delta$ on the complement (which is the closed set $K$ and contains the end $((0,1)_{r'}\times Y, r'\alpha_0)$) for $\delta$ small. We have five classes of Hamiltonian orbits similarly as before.

Now we choose an almost complex structure $J$, compatible with the symplectic form, cylindrical on the positive end of $\widehat{W}$, and cylindrical on the negative end of $(0,1)_{r'}\times Y$, which prescribes negative asymptotics of curves with finite Hofer energy have to be Reeb orbits of $\alpha_0$. Let $J_0$ denote the cylindrical almost complex structure, but defined on the whole symplectization $(0,\infty)_{r'}\times Y$. Let $H\in \cH(W)$, and $x_+,x_-$ are two non-degenerate Hamiltonian orbits of $H$. We define the compactified moduli spaces for an ordered multiset $\Gamma=\{\gamma_1,\ldots,\gamma_m\}$ of Reeb orbits of $(Y,\alpha_0)$.
\begin{equation}\label{eqn:partial}
\cM_{W,H}(x_+,x_-;\Gamma)=\overline{\left\{
\begin{array}{c}
     u:\R_s\times S^1_t \backslash \{p_1,\ldots,p_m\}\to \widehat{W} \\
     \text{asymptotic markers at } p_i
\end{array}
\left|
\begin{array}{c}
    \partial_s u+J(\partial_t u-X_{H})=0, \\
     \displaystyle \lim_{s\to \infty} u(s,\cdot)=x_+,\\
     \displaystyle \lim_{s\to -\infty} u(s,\cdot)=x_-,\\
     u(p_i)=\gamma_i.
\end{array}
\right.\right\}/\R},
\end{equation} 
Here, the $\R$-translation on the domain moves the punctures along with asymptotic markers. 

\begin{figure}[htb]
\begin{center}
{\small
\begin{overpic}[scale=0.4]{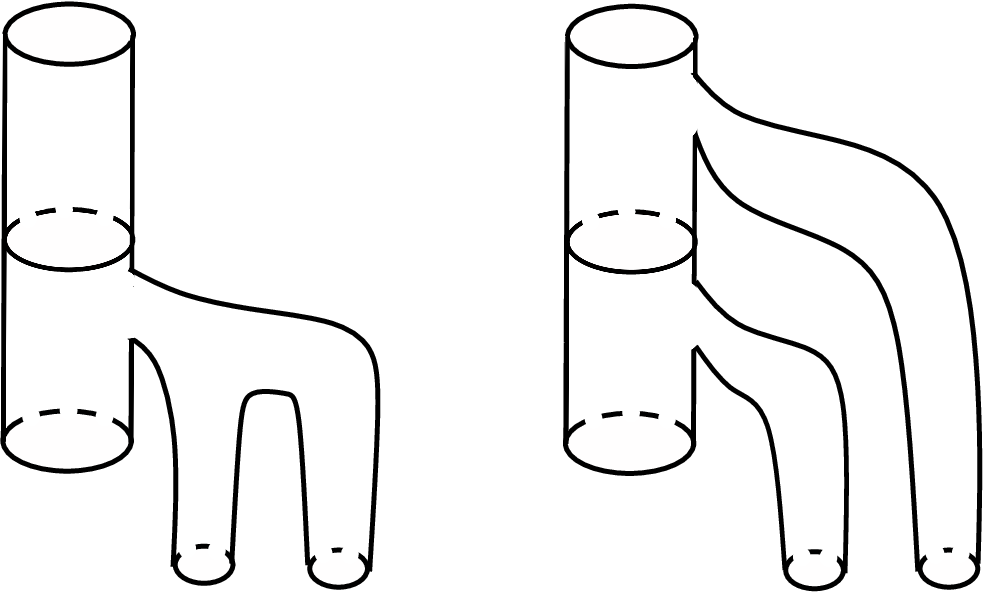}
\put (14, 58) {$x_+$}
\put (14,35) {$x_0$}
\put (10,10) {$x_-$}
\put (20,-4) {$\gamma_1$}
\put (35,-4) {$\gamma_2$}
\put (72,58) {$x_+$}
\put (72,35) {$x_0$}
\put (69,10) {$x_-$}
\put (82,-4) {$\gamma_1$}
\put (95,-4) {$\gamma_2$}
\end{overpic}
\hspace{1cm}
\begin{overpic}[scale=0.4]{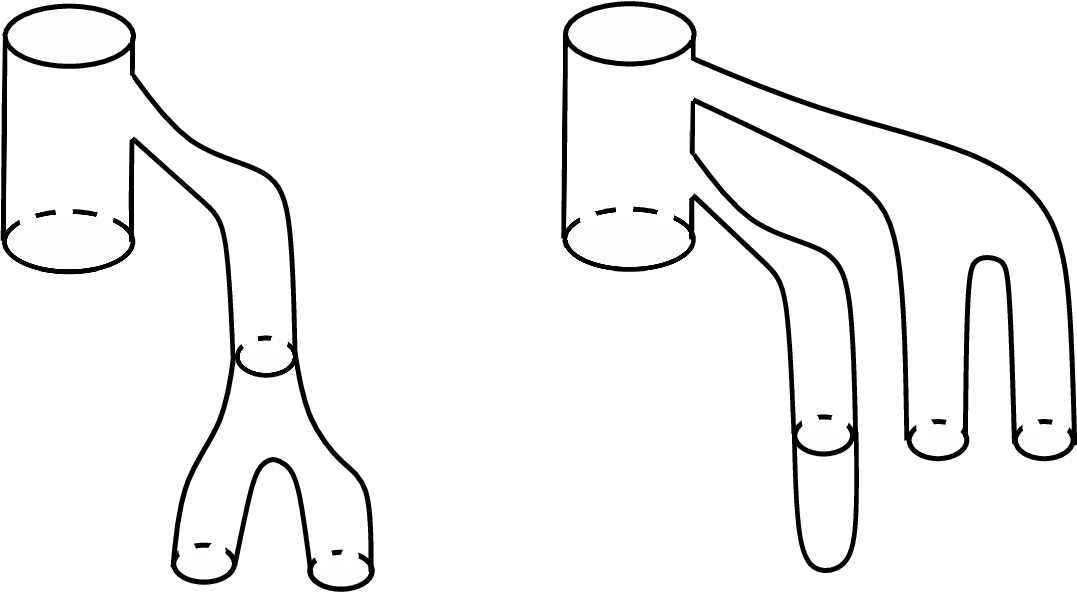}
\put (-6, 50) {$x_+$}
\put (-6, 30) {$x_-$}
\put (30,20) {$\gamma_0$}
\put (20,-4) {$\gamma_1$}
\put (32,-4) {$\gamma_2$}
\put (45,50) {$x_+$}
\put (45,30) {$x_-$}
\put (67,10) {$\gamma_0$}
\put (88,10) {$\gamma_1$}
\put (97,10) {$\gamma_2$}
\end{overpic}}
\caption{Codimension one degeneration, $x_*$ stands for Hamiltonian orbits and $\gamma_*$ are Reeb orbits}
\label{fig:deg}
\end{center}
\end{figure}
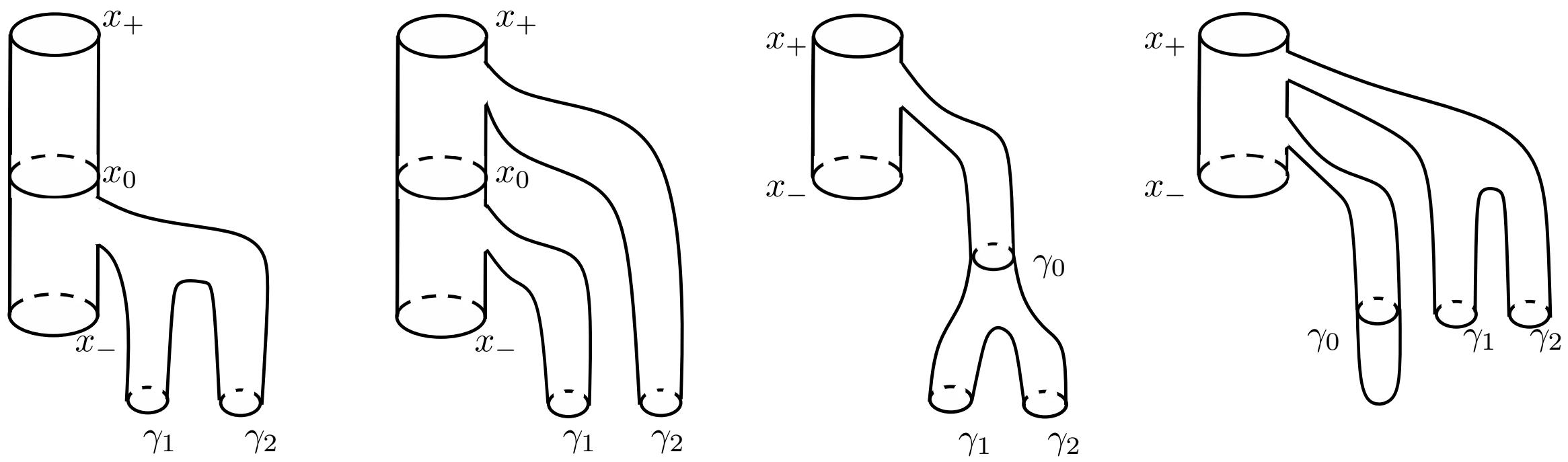

Now, we discuss the compactification and clarify some misunderstood part in the literature. For curves in \eqref{eqn:partial}, i.e.\ curves for defining the differential in $\RFH$, the underlying conformal structure is a parametrized cylinder $\R_s \times \bbS_t^1$ (to make sense of the Hamiltonian-Floer equation) with free interior marked points modulo the simultaneous $\R$-translation. The compactification involves the following  "codimension one" degeneration:
\begin{enumerate}
    \item The degeneration of the underlying conformal structures, namely:
    \begin{itemize}
        \item Some (but not all) of the marked points(i.e.\ negative punctures) can escape to $\infty$, leading to a Floer-type breaking involving two cylinders (modulo $\R$-translations on each cylinder) with marked points; See the second picture of Figure \ref{fig:deg}.
        \item Some of the marked points can collide, yielding an SFT curve in $\widehat{Y}$ as a ``sphere bubble'' with additional negative punctures. See the third picture of Figure \ref{fig:deg}.
    \end{itemize}
    \item Floer breakings at the two ends of $\R\times \bbS^1$, yielding a curve in \eqref{eqn:partial} with all the marked points and a curve in \eqref{eqn:partial} with no marked points. See the first picture of Figure \ref{fig:deg}.
    \item ``sphere bubbles" corresponding to the development of holomorphic planes in $\widehat{Y}$, which does not change the conformal structure on the main component. See the fourth picture of Figure \ref{fig:deg}.
\end{enumerate}

There is one more key point: a priori, when we solve the Hamiltonian-Floer equation on $\widehat{W}$ for the Hamiltonian $H$ with finite energy, the two ends can be asymptotic to both Hamiltonian orbits of $H$ and Reeb orbits of $(Y,\alpha_0)$. It is clear that we want to focus on non-constant Hamiltonian orbits of $H$. Note that both constant orbits and asymptotic Reeb orbits on $(Y,\alpha_0)$ has symplectic action $-H(K)=-\mu\delta$ (the $\lambda$ term is converging to $0$ for the asymptotic Reeb orbit case), which we may assume is sufficiently small. Now, once we fix an action window $(a,b)$, for $H(K)\gg 0$, i.e. $\mu\gg 0$, the orbits of types $F$ and $F_-$ can be thrown out, as there is no Floer breaking at those orbits for symplectic action reasons if we start with positive and negative asymptotic Hamiltonian orbits $x_+,x_-$ with symplectic action in $(a,b)$. Similarly, we can not have Floer-type breaking at the Reeb orbits, e.g.\ \cite[Figure 3]{Uebele} can not appear\footnote{In \cite{Uebele}, this was argued by the asymptotic behavior lemma \cite[Lemma 2.3]{CieliebakOancea2018} and the maximum principle \cite[Lemma 2.2]{CieliebakOancea2018} instead.}. See Figure \ref{fig:negative_energy} below, where the ``$ \Gamma$''-shaped curve can not exist for symplectic action reasons.

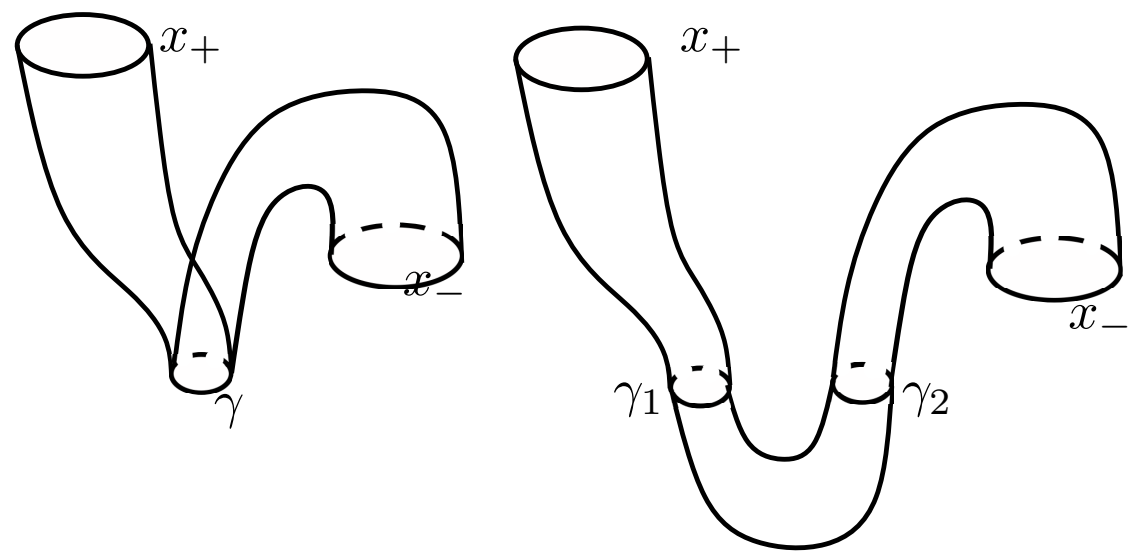
\begin{figure}[htb]
\begin{center}
{\small
\begin{overpic}[scale=0.4]{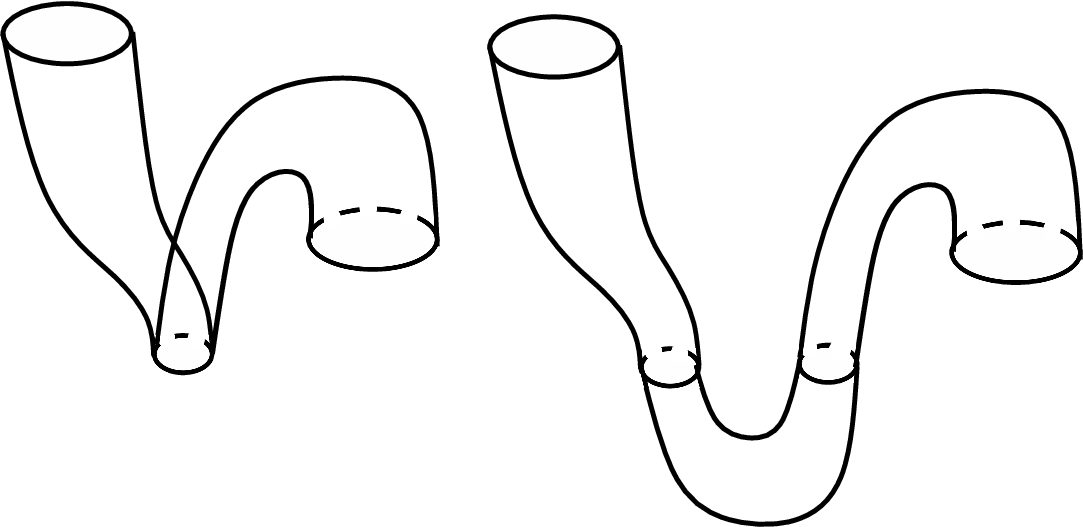}
\put (13, 45) {$x_+$}
\put (35,23) {$x_-$}
\put (18,12) {$\gamma$}
\put (60,45) {$x_+$}
\put (95,20) {$x_-$}
\put (54,13) {$\gamma_1$}
\put (80,13) {$\gamma_2$}
\end{overpic}}
\caption{Floer breaking involving a component with positive end asymptotic to a Reeb orbit and negative end asymptotic to a Hamiltonian orbit.}
\label{fig:negative_energy}
\end{center}
\end{figure}
As a consequence, the compactification $\cM_{W,H}(x_+,x_-;\Gamma)$ involves Floer breakings at Hamiltonian orbits (i.e.\ breakings at a Hamiltonian orbit, involving both a Floer differential breaking and separations of marked points in opposite directions of $\R$) of types $I_\pm,I_0$, as well as SFT breakings from developing holomorphic planes or SFT building trees from collision of marked points, each curve in $\widehat{Y}$ has exactly one positive puncture. Those curves breaking only once at a Hamiltonian orbit or a Reeb orbit form the expected boundary stratum of codimension $1$.

Now we fix an auxiliary choice $\Theta_0$ to define the contact DGA $\CC_*(Y,\Theta_0)$ of $(Y,\alpha_0,J_0)$, combining the virtual setup of Hamiltonian-Floer homology and contact homology by the works of Pardon \cite{Pardon:VFC,Pardon:Contact}, we then can make sense of the virtual counting $\#\cM_{W,H}(x_+,x_-;\Gamma)$ after choosing some auxiliary data $\Theta$ compatible with $\Theta_0$. Now we fix an augmentation $\epsilon$ of $\CC_*(Y,\Theta_0)$, for $a,b\notin \Spec(W,\lambda)$ and $\mu\gg 0$, we define, similarly as before, 
 $$\FC_*^{(a,b)}(H;\epsilon)=\displaystyle \bigoplus_{\substack{x\in \cP(H)\\ a\le \cA_H(x)<b}}\K \cdot x$$ with differential
 $$\partial(x_+)=\sum \frac{1}{\mu_{\Gamma}\kappa_{\Gamma}}\#\cM_{W,H}(x_+,x_-;\Gamma)\epsilon(q^{\Gamma})\cdot x_-.$$
It is a differential from the codimension $1$ boundary strata of the compactification, similar to the linearized contact homology, i.e.\ the breaking with holomorphic plane shall cancel with the degeneration from collision of marked points when weighted by the augmentation.

In \S \ref{ss:property} below, we will need TQFT properties of augmented $\SH$ as the filled $\SH$. Those properties, in particular, the various algebraic structures on $\SH_*^{\heartsuit}(W;F)$ are defined through counting Hamiltonian-Floer solutions over certain families of Floer data (conformal structures, lower-order perturbations etc.), typically a Floer data over a fixed Riemann surface with fixed cylindrical ends and Hamiltonian perturbations, e.g. the pair-of-pants product. For $\SH_*^{\heartsuit}(W;\epsilon)$, we count the same Hamiltonian-Floer equation with additional free marked points on the surface that can asymptote to Reeb orbits on $(Y,\alpha_0)$, then we have degeneration, in addition to those described by the Floer data and Floer differential breaking, (1) marked points colliding to one of the cylindrical ends forcing Floer breaking, (2) collision of marked points making a ``sphere bubble tree" from SFT buildings on $\widehat{Y}$ that is attached to the main component at one puncture, (3) developing holomorphic planes in $\widehat{Y}$. Similar, we before, when we work on some fixed action windows, for $H(K)\gg 0$, all Floer type breaking can only happen at Hamiltonian orbits of types $I^{\pm},I^0$. Now, the last two degeneration will cancel with each other once augmentation is added, and the first one shall contribute the differential with augmentation composed with certain structural map, the degeneration of Floer-data may contribute to other composition of structural maps as usual. To sum up,  for a fixed action window and $\mu \gg 0$, in addition to those Floer-type breaking in the filled case, we may only observe degeneration involving "bubble-trees" from SFT curves in $\widehat{Y}$. To have other degeneration involving curves with multiple positive punctures in $\widehat{Y}$ as rational SFT, we need to have the underlying degeneration of the conformal structure considered in the relevant family of Floer data. 
\begin{remark}
    By the discussion above, for fixed action windows and $H(K)\gg 0$ (i.e.\ $\mu\gg 0$), configuration in \cite[Figure 4]{Uebele} can not appear, as the upper-right component of the last figure in \cite[Figure 4]{Uebele} shall have negative energy by symplectic action reasons\footnote{$c_1$ in \cite[Figure 4]{Uebele} is a positive end and should be understood as a Floer-type breaking, it has symplectic action $-H(K)\ll 0$.}. As a consequence, instead of requiring product-index-positivity, i.e.\ $\mu_{\CZ}(\gamma)>3$ for all Reeb orbits as in \cite[Definition 3.15]{Uebele}, the well-definedness of product on a cobordism, i.e.\ \cite[Corollary 3.16]{Uebele}, needs only the weaker index positivity.
\end{remark}

Now, given a monotone homotopy of admissible Hamiltonians $H_s$, by counting the compactified moduli spaces:
\begin{equation}
\cM_{W,H_s}(x_+,x_-;\Gamma)=\overline{\left\{
\begin{array}{c}
     u:\R_s\times S^1_t \backslash \{p_1,\ldots,p_m\}\to \widehat{W} \\
     \text{asymptotic markers at } p_i
\end{array}
\left|
\begin{array}{c}
    \partial_s u+J(\partial_t u-X_{H_s})=0, \\
     \displaystyle \lim_{s\to \infty} u(s,\cdot)=x_+,\\
     \displaystyle \lim_{s\to -\infty} u(s,\cdot)=x_-,\\
     u(p_i)=\gamma_i.
\end{array}
\right.\right\}},
\end{equation} 
weighted by the augmentation $\epsilon$ of $\CC_*(Y,\Theta_0)$. We obtain a continuation map 
$$\FC_*^{(a,b)}(H_1;\epsilon)\to \FC_*^{(a,b)}(H_2;\epsilon).$$
If we fix $\alpha_0,J_0,\Theta_0,\epsilon$, by the same homotopy argument in Hamiltonian-Floer homology, the induced map on homology is a direct system on $H\in \cH(W)$. Therefore we define 
$$\SH_*^{(a,b)}(W;\epsilon):=\varinjlim_{H\in\cH(W)} \FH_*^{(a,b)}(H;\epsilon).$$
Similarly, we have the six versions $\SH_*^{\heartsuit}(W;\epsilon)$.

\subsubsection{Properties of Symplectic homology of augmented cobordism}\label{ss:property}

The followings are analogues \Cref{prop:basic_property}.
\begin{proposition}\label{prop:inv}
For a fixed $\epsilon$ (hence fixed $\alpha_0,J_0$ and $\Theta_0$), $\SH^{\heartsuit}_*(W;\epsilon)$ is an invariant, independent of the other choices, i.e.\ the Liouville form on $W$ up to Liouville homotopy and the almost complex structure $J$, as well as other auxiliary data used in the construction extending $\Theta_0$.
\end{proposition}
The proof is identical to the proofs of the filled version in \cite{CieliebakOancea2018} by carrying the weights from the augmentation $\epsilon$ through the proof, using the virtual setups in \cite{Pardon:VFC,Pardon:Contact}. More generally, one should expect invariance w.r.t.\ the data $\epsilon$ (hence $\alpha_0,J_0,\Theta_0$) up to homotopy. However, this requires establishing the DGA homotopy property for contact DGA (instead of chain complex homotopy established in \cite{Bao-Honda,Pardon:Contact}), which is yet missing in the construction of contact DGA. The issue is completely analogous to the status of linearized contact homology. As long as we do not touch the invariance w.r.t.\ to the contact DGA of the negative end, all the properties of $\SH_*^{\heartsuit}(W;\epsilon)$ can be established as the filled version in \cite{CieliebakOancea2018}, e.g.\ the following proposition. 

For the purpose of this paper, having an algebraic gadget with desired properties is sufficient, whether it is invariant for $\epsilon$ up to homotopy does not enter the argument. 
\begin{proposition}\label{prop:same}
\eqref{bp2} and \eqref{bp3} of \Cref{prop:basic_property}, \Cref{prop:dual,prop:product,prop:kill} and \Cref{cor:surj} holds for $\SH^{\heartsuit}_*(W;\epsilon)$.
\end{proposition}

We use $\RFH^{\heartsuit}_*(Y;\epsilon;\K)$ to denote the versions of the Rabinowitz Floer homology with coefficient field $\K$. This is defined on any trivial cobordism on $Y$, i.e.\ $\alpha_{\pm}$ could be differed by not a constant. However, we need a fix a contact form $\alpha_0$ and hence an almost complex $J_0$ on the negative end of $\widehat{Y}$ parametrized using the contact form $\alpha_0$, as well as auxiliary data $\Theta_0$. Those data are implicitly fixed in the augmentation $\epsilon$.

$\SH^{>0}(W;\epsilon;\K)$ and $\SH^{(a,b)}(W;\epsilon;\K)$ for $a\ge 0$ can be defined using Hamiltonians that are zero on $W\cup \partial_-W \times (0,1)$, e.g.\ \cite[]{AOT}. Let $H$ be a Hamiltonian on $\widehat{W}$ that is zero on $W\cup \partial_-W \times (0,1)$ and has slope $b$ on $\partial_+W\times (1+\epsilon,\infty)$ and $h$ is a Morse function on $W$ such that $\nabla_g h$ points out along $\partial W$ for a generic metric $g$ on $W$ such that $(h,g)$ is Morse-Smale. The connecting map 
$$\SH^{(0,b)}_*(W;\epsilon;\K)\to \SH^{=0}_{*-1}(W;\epsilon;\K)=H^{n+1-*}(W;\K)$$
can be defined by counting 

\begin{equation}\label{eqn:connecting}
\cM_{W,H}(x_+,z;\Gamma)=\overline{\left\{
\begin{array}{c}
     u:\R_s\times S^1_t \backslash \{p_1,\ldots,p_m\}\to \widehat{W} \\
     \text{asymptotic markers at } p_i \\
     \gamma:[0,\infty) \to W
\end{array}
\left|
\begin{array}{c}
    \partial_s u+J(\partial_t u-X_{H})=0, \\
     \gamma'(s) = \nabla_g h,\\
     \displaystyle \lim_{s\to \infty} u(s,\cdot)=x_+,\\
     \displaystyle \lim_{s\to -\infty} u(s,\cdot)=y\in W^{\circ},\\
     \displaystyle \lim_{s \to \infty} \gamma(s) = z,\\
      \gamma(0)=y, u(p_i)=\gamma_i.
\end{array}
\right.\right\}/\R},
\end{equation} 
weighted by the augmentation $\epsilon$, where $x_+$ is a non-constant Hamiltonian orbit of $H$ and $z$ is a critical point of $h$.

The Viterbo transfer map for $\SH$ was discussed in \cite[\S 5]{CieliebakOancea2018}. Here we just state a special case:
\begin{proposition}[{\cite[\S 5.1]{CieliebakOancea2018}}]\label{prop:Viterbo}
    Let $W$ be an exact cobordism and $W_0\subset W$ be an exact sub-domain, i.e.\ $W_0$ is an exact filling of $\partial W_0$. Let $\epsilon$ be an augmentation to $\partial_-W$, then we have Viterbo transfer map $\Phi$ making the following diagram commute:
    \begin{equation}\label{eqn:Viterbo}
        \xymatrix{
        \SH^{(0,b)}_*(W;\epsilon) \ar[r]^{\Phi}\ar[d] & \SH^{(0,b)}_*(W_0)\ar[d]\\
        H^{n+1-*}(W) \ar[r] & H^{n+1-*}(W_0)
        }
    \end{equation}
\end{proposition}

\subsubsection{Choices of coefficient}
There are several different choices of coefficient needed in this paper. For symplectic homology of cobordism $SH_*(W;F;\K)$ in \cite{CieliebakOancea2018}, any field coefficient is allowed as we only encounter moduli spaces without any orbifold points. For the symplectic homology of augmented cobordism $\SH_*(W;\epsilon;\K)$, as it is built on \cite{Pardon:Contact}, we must require that $\mathrm{char}(\K)=0$. In addition to those choices, we also require the following coefficients, which require modifications of the differentials in the definitions of $\SH_*(W;\epsilon;\K)$.

\paragraph{\textbf{Novikov coefficient}}
Let 
$$\Lambda_{\K} = \left\{\sum a_i T^{\lambda_i}\left| a_i\in \K, \lambda_i\in \R, \displaystyle \lim \lambda_i\to +\infty\right.\right\}$$
be the Novikov field over $\K$ with $\mathrm{char}(\K)=0$. We first define $\CC_*(\partial_-W,\alpha_0)_{\Lambda_{\K}}$ over $\Lambda_{\K}$ by
$$\partial_k(q_{\gamma}) = \sum_{[\Gamma],|\Gamma|=k} \frac{T^{\int \gamma^*\alpha_0-\sum_{\eta \in \Gamma}\int \eta^*\alpha_0}}{\mu_{\Gamma}\kappa_{\Gamma}} \# \cM_{\gamma,\Gamma}q^{\Gamma}.$$
Here $\CC_*(\partial_-W,\alpha_0)_{\Lambda_{\K}}$ is simply $\CC_*(\partial_-W,\alpha_0)\otimes_{\K} \Lambda_{\K}$, without taking one further completion from the $T$-filtration as in \cite[\S 2.6]{moreno2024rsft}. For example, elements like $\sum_{i=1}^\infty T^i q_{\gamma_i}$, for $\gamma_i\ne \gamma_j$ when $i\ne j$, do not lie in $\CC_*(\partial_-W,\alpha_0)_{\Lambda_{\K}}$, but in the completion in \cite[\S 2.6]{moreno2024rsft}. The completed theory is needed to discuss functoriality and Maurer-Cartan deformation from general cobordisms, but can be avoided here. By \cite[Corollary 3.15]{RSFT}, a strong filling would yield an augmentation over $\Lambda_{\K}$, the uncompleted version here, by counting curves in the filling. 

Now given a $\Lambda_{\K}$ augmentation $\epsilon$ of $\CC_*(\partial_-W,\alpha_0)_{\Lambda_{\K}}$, we can define 
$$\FC_*^{(a,b)}(H;\epsilon)=\displaystyle \bigoplus_{\substack{x\in \cP(H)\\ a\le \cA_H(x)<b}}\Lambda_{\K} \cdot x$$
with differential 
\begin{equation}\label{eqn:differential_Lam}
  \partial(x_+)=\sum \frac{T^{\cA_H(x_+)-\cA_H(x_-)-\sum_{\eta\in \Gamma}\int \eta^*\alpha_0}}{\mu_{\Gamma}\kappa_{\Gamma}}\#\cM_{W,H}(x_+,x_-;\Gamma)\epsilon(q^{\Gamma})\cdot x_-.  
\end{equation}
\begin{proposition}
    The weight in \eqref{eqn:differential_Lam} is positive.
\end{proposition}
\begin{proof}
    The weight $\cA_H(x_+)-\cA_H(x_-)-\sum_{\eta\in \Gamma}\eta^*\alpha_0$ is precisely the integral 
    $$\int_{\R\times S^1\backslash\{p_1,\ldots,p_m\}} ||\partial_s u || = \int_{\R\times S^1\backslash\{p_1,\ldots,p_m\}} ||\partial_t u - X_H||$$
    where the pseudonorm $||\cdot||$ is $\rd \lambda(\cdot,J\cdot)$ on $\widehat{W}\backslash (0,1)_r\times Y$ and $\rd \alpha_0 (\cdot, J \cdot)$ on $(0,1)_r\times Y$. This follows from the Stokes' Theorem with the setup of having $H$ is constant on $(0,1)_r\times Y$. The integrand is zero if $u(s, \cdot)$ is the same Hamiltonian orbit $x_+$ for all $s\in \R$, which is then not counted in the $0$-dimensional moduli space. 
\end{proof}

All other maps needed in setting up properties, e.g.\ continuation maps, etc., need to be modified similarly to take the weight into consideration. 

\begin{remark}
	In the construction of Rabinowitz-Floer theory with Novikov coefficient, there are at least two perpectives: (1) Taking the direct limit and inverse limit of slopes/Reeb period thresholds for the chain complex over the Novikov coefficient, which is the perspective in this paper. (2) Taking the completion w.r.t\ actions (which completes the Novikov coefficient as well) of the (homotopy) (co)limits of Hamiltonians, which is the perspective in \cite{ven}. The first perspective yields a topological invariant and is easier to be adapted to the augmented case, while the second perspective is sensitive to the contact form, as illustrated in \cite{ven}.
\end{remark}

\paragraph{\textbf{Twisted coefficient}}\label{ss:twisted}
We can further refine the coefficient to record the homology classes of the curves in the differential, i.e.\ twisted coefficients in \cite[\S 2.1]{AT}. A specialization of such twisted coefficients is the local system on the loop space, e.g.\ those used in \cite{subflexible}.

Let $W$ be an exact cobordism. We fix a basis of $H_1(W;\Q)$ and chain representatives $x_1,\ldots,x_m$. Now for any loop $x$ in $W$ with $[x]=\sum a_i[x_i]$ in $H_1(W;\Q)$, we $O(x)$ be the set of equivalence classes of rational $2$-chains $u$, such that $\partial u=x-\sum a_i x_i$ and $u_1\simeq u_2$ iff $u_1-u_2=0$ in $H_2(W;\Q)$. Then $O(x)$ is naturally a $H_2(W;\Q)$-torsor. We can consider $O(x)\times_{\K[H_2(W);\Q]}K(\K[H_2(W;\Q)])$, where $K(\ast)$ is the field of fractions. This makes sense as $H_2(W;\Q)$ is torsion-free and abelian, which implies that $\K[H_2(W;\Q)]$ is an integral domain. Then we can add the Novikov parameter $T$ to $O(x)\times_{\K[H_2(W);\Q]}K(\K[H_2(W;\Q)])$ and take the completion to obtain $G(x)$, which is a $\Lambda_{K(\K[H_2(W;\Q)]}$ module. The same construction yields $\Lambda_{K(\K[H_2(W;\Q)/R]}$ for a subspace $R\subset H_2(W;\Q)$. Let $R_0$ be the kernel of $H_2(Y:=\partial_-W;\Q)\to H_2(W;\Q)$. One can define the contact homology DGA of $Y$ over $\Lambda_{K(\K[H_2(Y;\Q)/R_0]}$ by recording both contact energy, i.e.\ $\int u^*\rd \alpha_0$ for the contact form, and the homology classes modulo $R_0$. Note that $H_2(Y;\Q)/R_0\to H_2(W;\Q)$ is injective, we can define the contact homology over  as well. Let $\epsilon$ be an augmentation in $\Lambda_{K(\K[H_2(W;\Q)]}$. Similarly, by counting Floer cylinders in $W$ with negative punctures and recoding the energy as in \eqref{eqn:differential_Lam} as well as the homology class, we obtain the augmented $\SH$ using coefficient in $\Lambda_{K(\K[H_2(W;\Q)]}$.

\begin{remark}
    In the case of strong or even stable cobordisms, except for strong fillings, additional completions on chain complex instead of the coefficient are needed to reflect Gromov compactness and the non-exactness of the symplectic form near the boundary, see \cite[\S 2.1]{AT} or \cite[\S 3.2.3]{APT}.
\end{remark}

\section{Detecting the topology from Rabinowitz Floer homology}\label{s:top}
Let $X$ be the symplectic manifold in \Cref{thm:main}. $Y$ is the Boothby-Wang contact manifold constructed from $(X,\omega)$ with first Chern class an integral lift of $[-\omega]$ in $H^2(X;\Z)$. More precisely, there is a connection form $\mathbf{i}\alpha$ on the associated $U(1)$-bundle, i.e., the circle bundle, such that $\rd \alpha = 2\pi \Proj^*\omega$, where $\Proj$ is the projection of the bundle. We use $L_{\omega}$ to denote the associated line bundle. The form $\alpha$ endows $Y$ with the celebrated Boothby–Wang contact form, such that the Reeb vector field is $\partial_{\theta}$ and all simple Reeb orbits have period $2\pi$.

We use $\GW_{X,1}(\pt,\pt)$ to denote the genus $0$ Gromov-Witten invariants of degree $1$ curves, i.e.\ rational curves $u$ with $\int u^*\omega=1$, passing through two generic points. Namely, $\GW_{X,1}(\pt,\pt) = \sum_{A,\la A,\omega\ra=1} \GW_{X,A}(\pt,\pt)$.
We start the proof of \Cref{thm:main} with the following seemingly stronger assumption:
\begin{assumption}
We assume that $\GW_{X,1}(\pt,\pt)\ne 0$ or $\GW_{X,1}(\pt,\pt)=\pm 1$ for now in the conditions of \Cref{thm:main}.
\end{assumption} 
This assumption can free us from using the twisted coefficients in \S \ref{ss:twisted}. In fact, by \cite[Theorem 1.3]{AP25}, $H_2(X;\Q)$ is generated by $A$, the assumption above is not actually stronger. However, this reduction does not apply to Theorem \ref{thm:tangency} and the twisted coefficient is necessary.

\subsection{Lifting curves from $X$ }
The idea is lifting closed holomorphic spheres in $\GW_{X,1}(\pt,\pt)$ to punctured holomorphic curves in the symplectization $\widehat{Y}$, and then interpreting them as curves counted in the Rabinowitz Floer homology. 
Let $p,q$ be two points in $X,Y$ respectively; we choose an $\omega$-compatible almost complex structure $J_X$, such that curves contributing to  $\GW_{X,1}(p,\Proj(q))$ are cut out transversely. This is possible for generic $J_X$ as $\omega$ is an integral form and those curves have minimal symplectic area. We use $\gamma_p$ to denote the simple Reeb orbit on $(Y,\alpha)$ over the point $p$. We also fix a parametrization of $\gamma_p$ by $[0,2\pi]/0\sim 2\pi$ starting at a fixed base point $\overline{p}$. We use $\overline{\gamma}_p$ to denote this parametrized Reeb orbit. $J_X$ induces an almost complex structure $J_Y$ on $\widehat{Y}=\R_s\times Y$ by $J_Y(\partial_s)=R=\partial_{\theta}$ and $J_Y|_{\xi}=\Proj^*J_X$ as $\rD \Proj$ induces an isomorphism from $\xi_x$ to $TX|_{\Proj(x)}$.  We consider the following moduli space:
\begin{equation}\label{eqn:M_Y}
\cM_{\widehat{Y}}(\gamma_p,q):=\left\{u:\C \to \widehat{Y}\left| \overline{\partial}_{J_Y}u=0, u(0) = (0,q)\in \R\times Y, \lim_{s\to \infty} u(e^{s+\mathbf{i}t}) = \overline{\gamma}_p(t)  \right.\right\}/\R,    
\end{equation}
where $\R$ acts on $\C$ by dilations. We use $\cM_{\widehat{Y}}(\gamma_p,q;\Gamma)$ be the version of \eqref{eqn:M_Y} with additional negative punctures asymptotic to the Reeb orbits collection $\Gamma$. Because of the point constraint from $(0,q)$, there is no translation on the target and hence the SFT compactification behave like the SFT compactification for curves in a (trivial) cobordism. 

\begin{proposition}\label{prop:lift}
    For $p\ne \Proj(q)$, $\cM_{\widehat{Y}}(\gamma_p,q)$ is the SFT compactification of itself, cut out transversely, and $\#\cM_{\widehat{Y}}(\gamma_p,q) = \GW_{X,1}(p,\Proj(q))\ne 0$. If $\Gamma\ne \emptyset$, $\cM_{\widehat{Y}}(\gamma_p,q;\Gamma)$ is empty.
\end{proposition}
\begin{proof}
    To see  $\cM_Y(\gamma_p,q)$ is the SFT compactification of itself, $\gamma_p$ has minimal period so there is no SFT breaking from the top. If there are degeneration from the bottom, the curve in the cobordism level passing through $(0,q)$ must be a trivial cylinder as $\gamma_p$ has minimal period. However, since $p\ne \Proj(q)$, this is impossible. In particular,  for $\Gamma\ne \emptyset$, the assumption that $p\ne \Proj(q)$ implies that there is no curves in $\cM_Y(\gamma_p,q;\Gamma)$. 

    We may identify $L_{\omega}\backslash 0$ with $\widehat{Y}$ as an almost complex manifold with the almost complex structure on $L_{\omega}$ induced from the connection $\mathbf{i}\alpha$. A curve $u$ in $\cM_Y(\gamma_p,q)$ is equivalent to a $J_X$-holomorphic rational curve $v:\CP^1\to X$ with $\int v^*\omega=1$ such that $v(\infty)=p,v(0)=\Proj(q)$ and a meromorphic section $s$ of $v^*L_{\omega}$ such that $s(\infty)$ is a simple pole and $\displaystyle \lim_{x\in \R \to \infty}\frac{s(x)}{|s(x)|}=\overline{p}\in L_{\omega}|_{v(\infty)}$ and $s(0)=q\in L_{\omega}|_{v(0)}$, where the correspondence is from $v=\Proj(u)$. The direction condition at $\infty$ fixes the $S^1$-gauge when we lift a curve from $X$ to $Y$. In particular, we have a one-to-one correspondence between $\cM_Y(\gamma_p,q)$ and curves in $\GW_{X,1}(p,\Proj(q))$. The transversality for $\cM_Y(\gamma_p,q)$ follows from the transversality of $\GW_{X,1}(p,\Proj(q))$, as the linearized Cauchy-Riemann operator is triangular, whose diagonal consists of two surjective operators: the linearized Cauchy-Riemann operator for $\GW_{X,1}(p,\Proj(q))$ and $\overline{\partial}$ on $v^*L_{\omega}$ subjected to the pole condition. Also, as $\overline{\partial}$ on $v^*L_{\omega}$ is complex linear, we have that they have the same sign: $\#\cM_Y(\gamma_p,q) = \GW_{X,1}(p,\Proj(q))\ne 0$.
\end{proof}

We need to perturb the contact forms on $Y$ to produce a trivial cobordism, such that we can single out the minimal orbits on the positive boundary and the negative boundary is non-degenerate to meet the setup in \S \ref{s:RFH}. We will do this in two steps.

First, we choose a Morse function $f$ on $X$, with a unique local minima at $p$ with minimum value $1$. We assume $f-1$ is $C^2$-small. If we consider the contact form $f\alpha$ on $Y$, $\gamma_p$ becomes the unique Reeb orbit with minimal period $2\pi$ and is non-degenerate. We use $W_f$ to denote the cobordism from $(Y,\alpha)$ to $(Y,f\alpha)$, although strictly speaking, the cobordism becomes ``singular" along $\gamma_p$, the completion $\widehat{W}_f$ still make sense, which is just $\widehat{Y}$ but equipped with a different collar coordinate at the positive end. We require the complex structure $J$ on $\widehat{W}_f$ to have the property that $J$ is $J_Y$ near the negative end and is compatible with $f\alpha$ on the positive end, tamed by $\rd \lambda =\rd(r\alpha)$ on $W_f$. In particular, we have a $2$-form $\overline{\omega}$, which is $\rd \lambda$ on $W_f$ and is $\rd \alpha$ and $\rd f\alpha$ on two ends. Then we have $\int u^*\overline{\omega}\ge 0$ for a $J$-holomorphic curve $u$, which is precisely the difference of contact actions of asymptotic Reeb orbits at two ends. By putting $q$ on $\partial_-W_f$, we define the moduli spaces  $\cM_{W_f}(\gamma_p,q)$,  $\cM_{W_f}(\gamma_p,q;\Gamma)$ as before. 

\begin{proposition}\label{prop:perturb}
    There exists $J$ satisfying the above conditions, such that $\cM_{W_f}(\gamma_p,q)$ is the SFT compactification of itself, cut out transversely, and $\#\cM_{W_f}(\gamma_p,q) = \GW_{X,1}(p,\Proj(q))\ne 0$. If $\Gamma\ne \emptyset$, $\cM_{W_f}(\gamma_p,q;\Gamma)$ is empty.
\end{proposition}
\begin{proof}
    The compactness and emptiness follow from the same argument as in \Cref{prop:lift} as $\int u^*\overline{\omega}\ge 0$. To prove the non-vanishing result, we can apply neck-stretching at $\partial_-W_f$, then curves in $\cM_{W_f}(\gamma_p,q)$ break into curves in \Cref{prop:lift} and a trivial cylinder in $\widehat{W}_f$ connecting $\gamma_{p}$ on $\partial_{\pm}W_f$. Since the trivial cylinder is cut out transversely with algebraic count $1$, the claim follows. 
\end{proof}

Now we perturb $\alpha$ on $\partial_-W_f$ to $\alpha_0<\alpha$, such that $\alpha_0$ is non-degenerate. Let $W$ denote the trivial cobordism from $(Y,\alpha_0)$ to $(Y,f\alpha)$.
\begin{proposition}\label{prop:perturb2}
There exists a perturbed $\alpha_0$ and $J$, such that $\cM_W(\gamma_p,q)$ is the SFT compactification of itself, cut out transversely, and $\#\cM_W(\gamma_p,q) = \GW_{X,1}(p,\Proj(q))\ne 0$. $\cM_W(\gamma_p,q;\Gamma)$ is empty if $\Gamma \ne \emptyset$.
\end{proposition}
This follows from \Cref{prop:perturb}, for otherwise, we can have a sequence $\alpha_i$ and $J_i$ such that \Cref{prop:perturb2} fails and $(\alpha_i,J_i)$ converges to those in \Cref{prop:perturb}, which will violate the compactness and openness of transversality. 

Now we turn on a Hamiltonian $H$ on $\partial_+W=(Y,f\alpha)$. It is mostly simple to use an autonomous Hamiltonian for now. That is, we use a Hamiltonian $H$ to $\widetilde{W}$, that is zero on $W\cup \partial_-W \times (0,1)$ and is in the form of $H(r)$ on $\partial_+W\times (1,\infty)_r$ with $H'(r)=2\pi^+$ for $r>1+\delta$ and $H''(r)>0$ for $r\in (1,1+\delta)$. Here we use $a^+$ to denote a number $>a$ such that $a^+-a\ll 1$. Then there is a unique $S^1$-family of Hamiltonian orbits $\overline{\gamma}_p$, corresponding to $\gamma_p$, which is the only non-constant Hamiltonian orbit family. Following the cascades construction of the symplectic cohomology, e.g.\ \cite[\S 2.1]{RP} or \cite{BO09}, $\overline{\gamma}_p$ provides two generators $\check{\gamma}_p$ and $\hat{\gamma}_p$ to the Hamiltonian-Floer theory. Following the setup in \cite[\S 2.1]{RP} along with the possible negative punctures asymptotic to Reeb orbits, we define $\cM_{W,H}(\check{\gamma}_p,z;\Gamma)$ and $\cM_{W,H}(\hat{\gamma}_p,z;\Gamma)$ analogous to \eqref{eqn:connecting} for the critical point $z$ of a Morse function $h$ on $W$. We may assume $q\in W$ is the unique local minimum of $h$ on $W$.
Because $\overline{\gamma}_p$ is the only non-constant Hamiltonian orbit, we shall not see a non-trivial cascade as in \cite[Figure 1]{RP}. It is clear that $\SH^{(0,2\pi^+)}(W;\epsilon)$ for any $\epsilon$ is spanned by $\check{\gamma}_p,\hat{\gamma}_p$ with zero differential. 
\begin{proposition}\label{prop:Ham_curve}
    We may choose generic $J$ that is $S^1$-independent, such that 
    \begin{enumerate}
        \item\label{1} $\# \cM_{W,H}(\check{\gamma}_p,q)=\# \GW_{X,1}(p,\Proj(q))$ and  for the unique local minimum $q$ of $h$.
        \item\label{2} $\# \cM_{W,H}(\hat{\gamma}_p,z;\Gamma)= 0$, for any critical point $z$ of $h$ on $W$;
        \item\label{3} $\cM_{W,H}(\check{\gamma}_p,q;\Gamma) = \emptyset$ for $\Gamma\ne \emptyset$ and local minimum $q$.
    \end{enumerate}
\end{proposition}
\begin{proof}
    \eqref{1} follows from applying neck-stretching along $\partial_+W$ and \Cref{prop:perturb2} plus the trivial cylinder connecting Hamiltonian orbits and Reeb orbits in \cite[Proposition 4.7]{intersection}. \eqref{2} follows from $S^1$-equivariant transversality as in \cite[Proof of Proposition 2.9]{RP}. \eqref{3} follows from neck-stretching and  \Cref{prop:perturb2}. 
\end{proof}

\subsection{Rabinowitz Floer homology of $Y$}\label{ss:RFH_Y}
\begin{proposition}\label{prop:augmentation}
    The contact DGA of $\partial_-W$ for coefficient $\K=\Lambda$ or $\C$ admits an augmentation. 
\end{proposition}
\begin{proof}
In the former case, we can consider the natural strong filling $W$ of $Y$ from the total space of the line bundle. The counting of rigid holomorphic planes in $\widehat{W}$ weighted by a version of energy would give an augmentation over $\Lambda$, see \cite[Corollary 3.1]{RSFT}. In this case, the definition of $\CC_*(Y)$ and $\SH_*(W;\epsilon)$ requires weights from energy as well. Alternatively, the fact that $Y$ admits a strong filling implies that $\CH_*(Y)\ne 0$ over the field $\K=\C$. Then by \cite[Lemma 3.7]{AOT}, $\CC_*(Y)$ admits an ($\Z/2$-graded) augmentation $\epsilon$ over $\C$. 
\end{proof}

We first work with the case $\K=\C$ to avoid weights.

\begin{proposition}\label{prop:vanish}
    We have $\C^2=\RFH^{(0,2\pi^+)}_*(Y;\epsilon;\C):=\SH_*^{(0,2\pi^+)}(W;\epsilon;\C)$ and $\RFH^{\le 0}_*(Y;\epsilon;\C)\to \RFH^{<2\pi^+}_*(Y;\epsilon;\C)$ is zero.  In particular, $\RFH_*(Y;\epsilon;\C)=0$.
\end{proposition}
\begin{proof}
    By \Cref{prop:Ham_curve}, we have $\RFH^{(0,2\pi^+)}_*(Y;\epsilon;\C)  \to H^{n+2-*}(Y;\C)\to H^0(Y;\C)$ is surjective. By \Cref{cor:surj}, we have $\RFH^{(0,2\pi^+)}_*\to \RFH^{\le 0}_{*-1}(Y;\epsilon;\C)$ is surjective. Then by tautological exact sequence, $\RFH^{\le 0}_*(Y;\epsilon;\C)\to \RFH^{<2\pi^+}_*(Y;\epsilon;\C)$ is zero. In particular, the unital ring map $\RFH^{\le 0}_*(Y;\epsilon;\C)\to \RFH_*(Y;\epsilon;\C)$ is zero and hence $\RFH_*(Y;\epsilon;\C)=0$.
\end{proof}

\begin{proposition}\label{prop:kernel}
    $\RFH^{(0,2\pi^+)}_*(Y;\epsilon;\C)\to H^{n+2-*}(Y;\C)$ has a kernel of dimension $1$.
\end{proposition}
\begin{proof}
    This follows from \eqref{2} of \Cref{prop:Ham_curve}.
\end{proof}

\begin{proposition}\label{prop:rational_sphere}
    $\dim_{\C} H^*(Y;\C) =2$.
\end{proposition}
\begin{proof}
By \Cref{prop:vanish}, $\RFH^{(0,2\pi^+)}_*(Y;\epsilon;\C) \to \RFH^{\le 0}_{*-1}(Y;\epsilon;\C)$ is surjective and hence we have $\dim_{\C}\RFH^{\le 0}_{*}(Y;\epsilon;\C) \le 2$. The vanishing of $\RFH_*(Y;\epsilon;\C)$ implies that $\dim_{\C}\RFH^{>0}_{*}(Y;\epsilon;\C) \le 2$. Moreover, by \Cref{prop:basic_property,prop:dual}, we must have $1\le \dim_{\C}\RFH^{>0}_*(Y;\epsilon;\C)=\dim_{\C}\RFH^{<0}_*(Y;\epsilon;\C)\le 2$. If $1=\dim_{\C}\RFH^{>0}_*(Y;\epsilon;\C)=\dim_{\C}\RFH^{<0}_*(Y;\epsilon;\C)$, then by \Cref{prop:basic_property,prop:vanish}, we must have $\dim_{\C} H^*(Y;\C)=2$. If  $2=\dim_{\C}\RFH^{>0}_*(Y;\epsilon;\C)=\dim_{\C}\RFH^{<0}_*(Y;\epsilon;\C)$, we must have $2=\dim_{\C}\RFH^{\le 0}(Y;\epsilon;\C)$ by the vanishing of $\RFH_*(Y;\epsilon;\C)$ and the surjection $\C^2=\RFH^{(0,2\pi^+)}_*(Y;\epsilon;\C) \to \RFH^{\le 0}_{*-1}(Y;\epsilon;\C)$ is necessarily an isomorphism. As a consequence, $\RFH^{(0,2\pi^+)}_*(Y;\epsilon;\C) \to \RFH^{>0}_*(Y;\epsilon;\C)$ is an isomorphism. Hence by \Cref{prop:kernel}, the map $\RFH^{>0}_*(Y;\epsilon;\C)\to H^{n+2-*}(Y;\C)$ has a kernel of rank $1$. Then from the spectral sequence and duality in \Cref{prop:basic_property,prop:dual}, the vanishing of $\RFH_*(Y;\epsilon;\C)$ implies that $\dim_{\C}H^*(Y;\C)=2$ and there is a differential from $\ker \delta$ to $\coker \delta^{\vee}$ in the spectral sequence in \Cref{prop:basic_property} in order to have vanishing $\RFH$ in \Cref{prop:vanish}.
\end{proof}

\subsection{General $\K$}

Now that $Y$ is a rational homology sphere by \Cref{prop:rational_sphere}, hence $X$ is a rational homology $\CP^n$ by the Gysin exact sequence. By our assumptions in \Cref{thm:main}, we have $c_1(A)=n+1$ for any curve class $A$ contributing to $\GW_{X,1}(\pt,\pt)\ne 0$, and hence $c_1(X)=(n+1)\omega$ in $H^2(X;\Q)$. With such a large monotonicity constant, for the  Boothby-Wang contact form (which is Morse-Bott) on $Y$, the lowest SFT degree is $2n>1$, see e.g.\ \cite[Proposition 3.3]{LZ26}. Moreover, the first Chern class of $(Y,\xi)$ is rationally zero and $H_1(Y;\Z)$ is torsion, as a consequence, all Hamiltonian orbits on $\widehat{Y}$ has a well-defined $\Q$-valued Conley-Zehnder index, see e.g.\ \cite[Proposition 3.8]{FZ}.

\begin{proposition}\label{prop:no_negative}
	For all moduli spaces in $\widehat{W}$ without negative punctures whose virtual dimension is at most $1$, we may assume the compactification does not involve curves with negative punctures.  In particular, those moduli spaces are cut out transversely for a generic choice of almost complex structure and are manifolds.
\end{proposition}
\begin{proof}
    If not, we shall have curves with negative punctures asymptotic to the Morse-Bott family of Reeb orbits. Since the lower-semicontinuous SFT degree of those families is at least $2n>1$. At least one component is the limit must have a negative virtual dimension. As the curves enjoy the somewhere injective property with the Hamiltonian perturbation, a generic choice of $J$ ensures such curves do not appear in the compactification. As a consequence, for a generic choice $J$, the original moduli space is compact and is a manifold.
\end{proof}

As a corollary of \Cref{prop:no_negative}, we can redo the definition of $\RFH^{\heartsuit}_*(Y;\epsilon)$ but using holomorphic curves in $\widehat{Y}$ without negative punctures, with all the properties discussed in \S \ref{s:RFH} and for any field coefficient $\K$. This is because moduli spaces of virtual dimension one without any negative punctures asymptotic to Reeb orbits have the Gromov-Floer compactification without involving SFT-type degeneration by \Cref{prop:no_negative}. Such invariants already appeared in \cite[\S 9.5]{CieliebakOancea2018}, which is the augmented theory with trivial augmentation. We use $\RFH^{\heartsuit}_*(Y;\K)$ to denote those invariants, they satisfy the same properties in \S \ref{ss:RFH_Y}.

\begin{proposition}\label{prop:GW}
    If $\GW_{X,1}(\pt,\pt)=\pm 1$, then $Y$ is an integer homology sphere.
\end{proposition}
\begin{proof}
   Similar to the proof of \Cref{prop:rational_sphere}, $\dim_{\K}H^*(Y;\K)=2$ for any field $\K$, in particular, for $\K=\Q$ and any finite field. Then by the universal coefficient theorem. $Y$ is an integer homology sphere. 
\end{proof}

\begin{proposition}\label{prop:twisted}
    \Cref{prop:rational_sphere,prop:GW} holds, when we replace $\GW_{X,1}(\pt,\pt)$ by $\GW_{X,A}(\pt,\pt)$ in the assumptions of \Cref{thm:main}.
\end{proposition}
\begin{proof}
By the Gysin exact sequence, we have $H_2(Y;\Q)\to H_2(X;\Q)$ induced by the inclusion of $Y$ into the unit disk bundle of the prequantization bundle is an injection and surjects onto $\ker \la \omega,\cdot \ra: H_2(X;\Q)\to \Q$. We can use the base change induced from the inclusion $H_2(Y;\Q)\to H_2(X;\Q)$ to discuss $\Lambda_{K(\K[H_2(X;\Q)])}$-coefficient $\RFH$, as the natural strong filling of the prequantization bundle produce a natural $\Lambda_{K(\K[H_2(X;\Q)])}$-augmentation. This allows us to differentiate the rational homology classes of the lifts of curves in $X$ whose homology class paired with $\omega$ is $1$. Then we can argue as in \Cref{prop:rational_sphere} first for $\K=\Q$ to conclude that $Y$ is a rational homology sphere and for $\K=\Z/p$ as in \Cref{prop:GW} to argue that $Y$ is an integer homology sphere. 
\end{proof}

\begin{remark}\label{rmk:alternative}
    One may consider moduli spaces in Proposition \ref{prop:lift}, but with the point constraint $q$ replaced by a loop, a submanifold, or a pseudo-cycle in $Y$. Then we can consider the universal family over such moduli spaces, which is a trivial $\CP^1$-bundle with two sections, since we can fix the $S^1$-gauge as in the proof of Proposition \ref{prop:lift}. The we can pull the pseudo-cycle to $\infty$ by a non-zero section of the universal family, whose image approaches the base point $\overline{p}$ on the Reeb orbit $\gamma_p$. Hence we can prove directly that $Y$ is a rational homology sphere. This is similar to the proof of \cite[Proposition 6.1]{CieZhou}, but on a contact manifold. The necessity and advantage for passing to $\RFH$ are reflected in the proof of Theorem \ref{thm:tangency} (see Remark \ref{rmk:fail}), which is the main force to derive filling rigidity results. 
\end{remark}

\subsection{Proof of \Cref{thm:main}}
By \cite[Theorem 1.3]{AP25}, as our condition in \Cref{thm:main} implies that $X$ is enumeratively rationally connected, we have that $\pi_1(X)$ is finite. We will prove that $X$ is actually simply connected. 

\begin{proposition}\label{prop:finite}
    Under either assumption of \Cref{thm:main}, we have $\pi_1(X)=0$
\end{proposition}
\begin{proof}
    By \cite[Theorem 1.3]{AP25}, $\pi_1(X)$ is finite. Let $\widetilde{X}$ be the universal cover, which is a compact symplectic manifold. Let $S_A$ denote the homotopy classes $f$ in $[S^2,X]$ such that $h(f)=A\in H_2(X;\Q)$ for the Hurewicz map $h$. We use $\GW_{X,f}(\pt,\pt)$ to denote the Gromov-Witten invariants of curves in homotopy class $f$. Then we have 
    $$\GW_{X,A}(\pt,\pt)=\sum_{f\in S_A} \GW_{X,f}(\pt,\pt).$$
    By lifting curves to $\widetilde{X}$, we have 
    \begin{equation}\label{eqn:cover}
        |\pi_1(X)|\cdot \sum_{\substack{g\in [S^2,\widetilde{X}] \\ \pi_*(g)=f}}\GW_{\widetilde{X},g}(\pt,\pt)=\GW_{X,f}(\pt,\pt),
    \end{equation}
     where $\pi_*:[S^2,\widetilde{X}]=\pi_2(\widetilde{X})=\pi_2(X)\to \pi_2(X)/\pi_1(X) = [S^2,X]$. Therefore we have 
     $$\sum_{g\in \widetilde{S}_A} \GW_{\widetilde{X},g}(\pt,\pt) \ne 0$$
    where $\widetilde{S}_A$ is the subset of $[S^2,\widetilde{X}]=\pi_2(\widetilde{X})=H_2(\widetilde{X};\Z)$ such that $\pi_*(g)=A\in H_2(X;\Q)$. We may partition the set $\widetilde{S}_A$ by their rational homology classes, i.e.\ neglecting the torsions. In particular, we can find a rational homology class $B$ such that 
    $\GW_{\widetilde{X},B}(\pt,\pt)\ne 0$. Therefore, $\widetilde{X}$ is compact and satisfies the assumptions of \Cref{thm:main}. In particular, $\widetilde{X}$ has the same rational cohomology as $\CP^n$ generated by the symplectic form. As the deck transformation fixes the symplectic class, the deck-transformation must have fixed points by the Lefschetz fixed point theorem. As a consequence, we must have $X$ is simply connected. 
\end{proof}

\begin{remark}
    Instead of using \cite[Theorem 1.3]{AP25}, one may prove $\pi_1(X)$ is finite directly. This would follow from the lifting argument for \eqref{eqn:cover} in the proof above, as \eqref{eqn:cover} would never be true if $|\pi_1(X)|$ is infinite. Note that even though $\widetilde{X}$ is non-compact, compactness for $\GW_{\widetilde{X},g}(\pt,\pt)$ still holds for the pullback of the almost complex structure from $X$. This, in some sense, is unwrapping the argument in \cite{AP25}.
\end{remark}

\begin{proof}[Proof of \Cref{thm:main}]
By the discussion above, $Y$ is a rational homology sphere by \Cref{prop:rational_sphere,prop:twisted}. By Gysin exact sequence, we have  $H^*(X;\Q)=\Q[x]/x^{n+1}$ as a ring for $x=[\omega]$.  Let $f$ be the map $X\to \CP^{\infty}=K(\Z,2)$ representing the generator $x\in H^*(X;\Z)$ in the free part. We may assume $f$ maps $X$ to the $2n$-skeleton of $\CP^{\infty}$, namely $\CP^n$, and still use $f$ to denote this map $X\to \CP^n$. Then $f$ induces an isomorphism on rational cohomology, as it is generated by $x$. Then $f$ is a rational homotopy equivalence, as $X,\CP^n$ are simply connected by \cite[Theorem 8.6]{zbMATH01510507}.If $\GW_A(\pt,\pt)=\pm 1$, then $Y$ is an integer homology sphere by \Cref{prop:GW,prop:twisted}. Therefore  $H^*(X;\Z)=\Z[x]/x^{n+1}$ as a ring by Gysin exact sequence. The same argument above shows that the classifying map for the generator $x\in H^2(X;\Z)$ represents a homotopy equivalence from $X$ to $\CP^n$. 

Since $A$ is the positive generator of $H_2(X;\Z)$ and $c_1(A)=n+1$, from that $\GW_{X,A}(\pt,\pt)\ne 0$, the quantum cohomology as a vector space over the Novikov ring is $\Q[x]/x^{n+1}\otimes \Lambda_t$ with $\deg(t)=2n+2$. Now $\GW_{X,A}(\pt,\pt)\ne 0$, implies $3$-point Gromov-Witten invariant $\GW_A(\pt,\pt,H)\ne 0$ by the divisor axiom, where $H\in H_{2n-2}(X;\Z)$ is an generator whose intersection with $A$ is $1$, i.e. $H$ is the Poincar\'e dual of generator $x$ of $H^2(X;\Z)$. $\int_X x^n =d>0 \in \Z$ implies that the Poincar\'e dual of a point is $\frac{1}{d}x^n$. We can interpret $\GW_A(\pt,H,\pt)\ne 0$ as $x^n \cup_Q x=d\GW_A(\pt,H,\pt) t $ in the small quantum group. This relation determines the quantum product by associativity. We have that the small quantum group of $X$ is the same as $\CP^n$ (after the change of coordinates $d\GW_{X,A}(\pt,\pt) t \mapsto t$). When $\GW_{X,A}=1$, we have $\int_X \omega^n = 1$, hence $H^*(X;\Q)\simeq H^*(\CP^n;\Q)$ as a Frobenius algebra. As $H^*(X;\Q)$ is generated by $H^2(X;\Q)$, we can apply the work of Kontsevich-Manin \cite[Theorem 3.1]{KM94} as well as Ruan-Tian \cite{RT95} to see that the big quantum cohomology/genus $0$ Gromov-Witten potential is also determined. 
\end{proof}

\begin{corollary}\label{cor:McDuff}
    Under the assumption of \Cref{thm:main}, if $X$ is $4$-dimensional, then $X$ is symplectomorphic to $(\CP^2,\omega_{\FS})$.
\end{corollary}
\begin{proof}
    Since $X$ is simply connected and rationally homotopy equivalent to $\CP^2$, by Freedman's result, $X$ must be homeomorphic to $\CP^2$. Moreover, from $\GW_{X,A}(\pt,\pt)\ne 0$ and $\la A,\omega \ra=1$, we have $A$ is the positive generator of $H_2(\CP^2;\Z)$ and $c_1(A)=3$. As a consequence, $X$ has no embedded symplectic sphere $S$ with self-intersection number $S\cdot S=-1$, i.e.\ $X$ is minimal. Since the topological intersection number $A\cdot A$ is $1$, by the curve adjunction formula  \cite[Theorem 2.51]{WenBook}, the curve must be an embedding. Then by \cite[Theorem 1.4]{zbMATH04192555}, $X$ is the standard symplectic $\CP^2$.
\end{proof}

\section{Proof of \Cref{thm:tangency}}\label{s4}
\subsection{$S^1$-equivariant $\SH$}
The $S^1$-equivariant version of Rabinowitz-Floer homology/symplectic homology of Liouville cobordism was defined in \cite[\S 8.2]{CieliebakOancea2018} for filled Liouville cobordisms.
The same construction can be generalized to augmented Liouville cobordisms as in \S \ref{ss:SH_aug}, by counting parametrized Floer cylinders \cite[\S 3.1]{dilation} with negative punctures weighted by the augmentation. One important feature of the $S^1$-equivariant version is the $u$-map coming from $H^*(BS^1)=\K[u]$ making the $\SH_{*}^{\heartsuit, S^1}(W,\epsilon)$ into a $\K[u]$-module.

Analogous to \eqref{eqn:connecting}, e.g.\ \cite[\S 5.2.2]{dilation}, by counting parametrized Floer cylinders with a point constraint for Hamiltonians that are zero on $W\cup \partial_-W \times (0,1)$ possibly with negative punctures weighted by the augmentation $\epsilon$, we get a map 
$$\delta_0: \SH^{(0,b),S^1}(W;\epsilon) \to \SH^{=0,S^1}(W)=H^*(W)\otimes \K[u,u^{-1}]/u\K[u] \to H^0(W)=\K.$$
By counting parametrized Floer equations with a point constraint with local tangency condition of order $k\ge 0$, we get a map 
$$\delta_k: \SH^{(0,b),S^1}(W;\epsilon)\to \K$$
supported in degree $n+1+2k$ for $2n=\dim W$. 

The Viterbo transfer map in \Cref{prop:Viterbo} has an $ S^1$-equivariant enhancement, which is compatible with $\delta_k$. This, in the filled case, is the usual Viterbo transfer map for the connecting map $\SH^{+,S^1}\to \SH^{=0,S^1}$ in \cite[\S 4.2]{Gutt}, with the connecting map related by $\delta_k$ by counting parametrized Floer cylinders with a tangency condition.
\begin{proposition}\label{prop:Viterbo_S1}
    Let $W_0$ be an exact subdomain in the exact cobordism $W$. 
    We have a commutative diagram 
    $$\xymatrix{
    \SH^{(0,b),S^1}_*(W;\epsilon)\ar[d]^{\delta_i} \ar[r]^{\Phi} & \SH^{+,S^1}_*(W_0) \ar[d]^{\delta_i} \\
    \K \ar[r]^{=} & \K}
    $$
    where the horizontal map $\Phi$ respects the $\K[u]$-structure.
\end{proposition}
In view of \cite{BO09}, $\SH_*^{(0,b),S^1}(W,\epsilon)$ is the (truncated) linearized contact homology of $\partial_+W$ w.r.t.\ the augmentation from the composition of the DGA map induced from the cobordism $W$ and $\epsilon$. After this identification, the transfer map along with the non-vanishing result for $\delta_i$ on a Darboux ball (by Proposition \ref{prop:tangency_u} below) is the analogues of those established by Siegel in \cite[\S 5.5]{Siegel}.

\subsection{Proof of \Cref{thm:tangency}}
For \Cref{thm:tangency}, similar to the proof of \Cref{thm:main}, we write $Y$ as the prequantization circle bundle over $X$ with first Chern class $-[\omega]$. We also assume that 
$$\GW_{X,1}(B,\pt,k)\ne 0$$
for all rational curves with area $1$. The general case for $\GW_{X,A}(B,\pt,k)\ne 0$ is completely similar to \Cref{thm:main} using the fraction field of the group ring.

\begin{proposition}\label{prop:tangency_u}
For any augmentation $\epsilon$ of $Y$,  the following holds:
\begin{enumerate}
    \item $\RFH_{*}^{(0,2\pi^+),S^1}(Y,\epsilon) \simeq H^{n-*}(X)$ as a $\K[u]$-module (in $\Z/2$ grading), where $u$ action on $H^{*}(X)$ is the cup product with $[\omega]$.
    \item The class $v\in \RFH_{*}^{(0,2\pi^+),S^1}(Y,\epsilon)$, which in $H^{*}(X)$ is Poincar\'e dual to $B$, such that $\delta_k(v) = \GW_{X,1}(B,\pt,k)$. 
\end{enumerate}
\end{proposition}
\begin{proof}
    This follows from a similar argument as in \cite[\S 5.2.1, 5.2.2]{dilation} using the Boothby-Wang Morse-Bott contact form on $Y$.
\end{proof}
\begin{remark}
    If $Y$ has vanishing rational first Chern class, i.e.\ $c_1(X)$ is proportional to $[\omega]$, then we have $\RFH_{*}^{(0,2\pi^+),S^1}(Y,\epsilon) \simeq H^{-*+n+2C}(X)$, where $C=c_1(X)/[\omega]$.
\end{remark}

\begin{example}\label{ex:sphere}
    A special case is when $X=(\CP^n,\omega_{\FS})$, then $Y=(S^{2n+1},\xi_{\std})$. Let $H_k$ be a complex dimension $k$ linear subspace in $\CP^n$ for $k<n$, then by \cite[Proposition 6.4]{CieZhou}, we have 
    $$\GW_{\CP^n,[H_1]}(H_k,\pt, k) =\GW_{\CP^n, 1}(H_k,\pt,k)=k!.$$
    Then, by Proposition \ref{prop:tangency_u}, we have 
    $$\RFH^{(0,2\pi^+),S^1}_*(S^{2n+1},\epsilon) = \la v_i \ra_{0\le i \le n}, \quad u\cdot v_i=v_{i-1} $$
    and 
    $$\delta_i(v_i)=i!,\quad i<n$$
    More generally, when $\epsilon$ is from the natural symplectic ball filling of $(S^{2n+1},\xi_{\std})$, we have 
    $$\RFH^{+,S^1}_*(S^{2n+1},\epsilon_{B^{2n+2}}) = \la v_i \ra_{0\le i \le \infty}, uv_i=v_{i-1}$$
    and 
    $$\delta_i(v_i)=i!, \quad \forall i<n.$$
    One can use $\GW_{\CP^{n+1},[H_1]}(H_n,\pt, n)=n!$ to conclude that $\delta_n(v_n)=n!$ as well.
\end{example}

\begin{proof}[Proof of \Cref{thm:tangency}]
    Let $\Phi$ be the Viterbo transfer map in \Cref{prop:Viterbo_S1} and $v\in \RFH_{*}^{(0,2\pi^+),S^1}(Y,\epsilon)$ from \Cref{prop:tangency_u}. \Cref{prop:Viterbo_S1} and non-vanishing results in \Cref{ex:sphere} imply that 
    $$\Phi(v) = \frac{\GW_{X,1}(B,\pt,k)}{k!} v_k +\sum_{j\ne k} a_j v_j.$$
    Let $\mu$ be the element in $\RFH_{*}^{(0,2\pi^+),S^1}(Y,\epsilon)$ representing the generator of $H^{2n}(X)$. Then $\la \omega^k,B\ra=m$ and \Cref{prop:tangency_u} implies that $u^{k}\cdot v = m\mu $. Then \Cref{prop:Viterbo_S1} and \Cref{ex:sphere} for $i=0$ implies that 
    $$\delta_0:\RFH_{*}^{(0,2\pi^+),S^1}(Y,\epsilon) \to \K, \quad \mu \mapsto \frac{\GW_{X,1}(B,\pt,k)}{m\cdot k!}.$$
    Applying \Cref{prop:tangency_u} reversely, we see that 
    $$\GW_{X,1}(\pt,\pt,0) = \frac{\GW_{X,1}(B,\pt,k)}{m\cdot k!}\ne 0.$$
    Then we reduce \Cref{thm:tangency} back to \Cref{thm:main}. The case for $\GW_{X,A}(B,\pt,k)\ne 0$ follows similarly with the group ring coefficient. 
\end{proof}

\begin{remark}\label{rmk:fail}
If one wishes to prove \Cref{thm:tangency} in a similar manner to \Cref{rmk:alternative}, one has to use the additional tangency condition in an essential way. This is because for quadratic hypersurface $Q^n\subset (\CP^{n+1},\omega_{\FS})$, we have a non-vanishing Gromov-Witten invariant $\GW_{Q^n,1}(S,\pt)$ without the tangency condition for $S\subset Q^{n}\subset \CP^{n+1}$ a complex line. A proof without Rabinowitz Floer homology as in \Cref{rmk:alternative} should explain this discrepancy effectively. More precisely, the homotopy to infinity in  \Cref{rmk:alternative} requires fixing the $S^1$-gauge, but the circle bundle over $S$ is non-trivial. In particular, there will be obstructions unless the dimension is very low, i.e.\ dimension $0,1$ and some dimension $2$ with boundary. Those low dimensional cases will be used in the proof of \Cref{thm:RP}. 
\end{remark}

\begin{proof}[Proof of Corollary \ref{cor}]
    By Theorems~\ref{thm:main} or~\ref{thm:tangency}, we have $H^2(X;\Q)=\Q$ and $\pi_1(X)=0$. Hence, $X$ has Picard number $1$. Then as $\la c_1(X),A\ra = (n+1) \la \omega,A \ra$, when $\omega$ is K\"ahler, we have $(X,\omega)$ is Fano with Fano index $n+1$ and $X$ is biholomorphic to $\CP^n$ by the Kobayashi-Ochiai Theorem. In the case of $\eta$ being K\"ahler, we may normalize $\eta$ to be primitive and integral. Then we have $[\eta]=\pm [\omega]$. As $\eta$ and $\omega$ give rise to the same almost complex manifold, in particular, the same orientation, we have $[\eta]^n = [\omega]^n$. Therefore when $n$ is odd, we have $[\eta]=[\omega]$ in cohomology and $(X,\eta)$ is Fano with Fano index $n+1$. The claim follows.
\end{proof}
When $\omega$ is K\"ahler, the case of Corollary~\ref{cor} arising from \Cref{thm:main} can alternatively be derived from Mori's cone theorem as follows: $X$ is projective by the integer assumption on $\omega$, $A$ is necessarily a curve class in an extremal ray, and the contraction of this ray must be a point because $\GW_{X,A}(\pt,\pt)\neq 0$\footnote{We implicitly assume that the symplectic Gromov-Witten invariant used in \Cref{thm:main} agrees with the algebraic Gromov--Witten invariant used here.}. 
Consequently, the Picard number of $X$ is $1$, so $X$ is Fano with Fano index $n+1$, and the result follows from Kobayashi-Ochiai. 
For Corollary~\ref{cor} under the assumptions of \Cref{thm:tangency} with $\omega$ K\"ahler, it is unclear whether one can still deduce that the Picard number is $1$ after the contraction in order to compute the Fano index as before. 
Theorems~\ref{thm:main} and~\ref{thm:tangency} bypass the cone theorem and gain some flexibility with respect to the symplectic form.	
\section{Applications to symplectic fillings}\label{s:filling}
\subsection{The Eliashberg-Floer-McDuff theorem}
We first give an alternative proof of the Eliashberg-Floer-McDuff theorem using \Cref{thm:main}.

\begin{corollary}[Eliashberg-Floer-McDuff]\label{cor:EMF}
    Let $W$ be an exact filling of $(S^{2n-1},\xi_{std})$ for $n\ge 3$, then $W$ is diffeomorphic to a $2n$-dimensional ball.  
\end{corollary}
\begin{proof}
    By the standard $h$-cobordism argument in the original proof \cite{zbMATH04186536}. It suffices to prove that $W$ is a homology ball and is simply connected. Let $B\subset (\CP^n,\omega_{FS})$ be a Darboux ball. Then $X=(\CP^n\backslash B) \cup W$ is an example satisfying all the assumptions in \Cref{thm:main}. Indeed, let $A\in H_2(\CP^{n-1};\Z)$ be the positive generator of the divisor $H=\CP^{n-1}$ at infinity,  we claim that $\GW_{X,A}(\pt,\pt) = 1$. Since the intersection number of $A$ with the divisor $H$ is $1$, by picking an almost complex structure on $X$ making $H$ holomorphic and two point-constraints $p,q\in H\subset X$, the moduli spaces for $\GW_{X,A}(p,q)$ and $\GW_{H,A}(p,q)$ are the same. By the analysis of the normal Cauchy-Riemann operator as in \cite[Lemma 2.3]{CieZhou}, we have $\GW_{X,A}(p,q)=\GW_{H,A}(p,q)=1$.  By \Cref{thm:main}, $X$ is homotopy equivalent to $\CP^n$, we can argue by the Mayer-Vietoris sequence and the van Kampen theorem that $W$ is homotopy equivalent to a point. 
\end{proof}

Let $\Theta_n$ be the group of homotopy $n$-spheres. By Cerf's theorem, for $n\ge 6$, we have $\Theta_n$ is the mapping class group $\pi_0(\mathrm{Diff}_+(S^{n-1}))$. Given a smooth manifold $M$, the inertia group $I(M)$ is the subgroup of $\Theta_n$ that consists of homotopy spheres $\Sigma$, such that $M\#\Sigma$ is diffeomorphic to $M$. We are interested in the quotient group $\Theta_{2n}/I(\CP^n)$. By \cite{zbMATH06828008}, we have $\Theta_{2n}/I(\CP^n)\ne 0$ for $n=4,5,7,8,9$.

\begin{proposition}
Let $\pi_0(\mathrm{Cont}(S^{2n-1},\xi_{std}))$ be the contact mapping class group. When $\Theta_{2n}/I(\CP^n)\ne 0$, then 
\begin{enumerate}
    \item either $\pi_0(\mathrm{Cont}(S^{2n-1},\xi_{std}))\to \pi_0(\mathrm{Diff}_+(S^{2n-1}))\to \Theta_{2n}/I(\CP^n)$ is not surjective,
    \item or there exists a symplectic, non-K\"ahler (in the sense that there is no K\"ahler structure) manifold $M$ that is homeomorphic to $\CP^n$.
\end{enumerate}
\end{proposition}
\begin{proof}
    If the map is surjective, there exists a non-trivial element in $\Theta_{2n}/I(\CP^n)$ represented by a contactomorphism $\phi$ of $(S^{2n-1},\xi_{std})$. Then, for a Darboux ball $B$,  $M:=(\CP^n\backslash B \cup_{\phi} B)$. By construction, $M$ is only homeomorphic to $\CP^n$. By a theorem of Hirzebruch, Kodaira and Yau, see \cite[Theorem 1.1]{TV}, if $M$ is K\"ahler, then $M$ is diffeomorphic to $\CP^n$. Therefore $M$ is a symplectic, non-K\"ahler manifold $M$ that is homeomorphic to $\CP^n$.
\end{proof}
Let $\phi$ be a diffeomorphism of $S^{2n-1}$, to check whether $\phi^*\xi_{std}$ is contact isotopic to $\xi_{std}$ seems to be a non-trivial task. It would be interesting to investigate the first example, namely the non-trivial element in $\Z/2=\pi_0(\mathrm{Diff}_+(S^{7}))=\Theta_8=\Theta_8/I(\CP^4)$. 

\subsection{Filling rigidities from Theorems \ref{thm:main} and \ref{thm:tangency}}
We explain the curves in hypersurfaces of low degrees. 
\begin{lemma}\label{lemma:GW}
    Let $S$ be a degree $k$ complex hypersurface in $\CP^n$ for $1\le k<n$ and $n\ge 3$, then $$\GW_{S,1}(B,p) = \la B, \omega^{k-1}\ra(k-1)!$$ for an oriented submanifold $B$ in $S$ of dimension $2(k-1)$. Let $W$ be a symplectic filling of $\partial N(S)$ where $N(S)$ is a tubular neighborhood of $S$, and $X$ is obtained from gluing $N(S)$ with $W$, then 
    $$\sum_A \GW_{X,A}(B,p,k-1)=\GW_{S,1}(B,p),$$
    where the sum is over the set of homology classes in $H_2(X,\Z)$ that are from $S$ and sent to the positive generator of $H_2(\CP^n;\Z)$.
\end{lemma}
\begin{proof}
    By \cite[Lemma 2.3]{CieZhou}, we have 
    $$\GW_{S,1}(B,p) = \GW_{\CP^n,1}(B,p,k-1)$$
    where the latter is given by $\la B,\omega^{k-1}\ra (k-1)!$ by \cite[Proposition 6.4]{CieZhou}. The second claim follows from \cite[Lemma 2.3]{CieZhou}, as the class $A$ has the property that $A\cdot S = k$ in $X$. 
\end{proof}

Let $N_D$ be a neighborhood of $D$ as a symplectic cap of $Y_D$. Let $W$ be a Liouville filling of $Y_{D}$ with a Liouville form $\mu$ such that $\mu=\lambda$ near $Y_D$, where $\lambda$ is the Liouville form on $\CP^n\backslash N_D\subset A_D$. Then $(X,\omega):=(N_{D}\cup_{Y_{D}} W, \omega_{\FS}|_{N_{D}}\cup \rd \mu)$ is a closed symplectic manifold.
\begin{proposition}\label{prop:integer}
    Under the assumption of \Cref{thm:homology_unique}, $(X,\omega)$ above satisfies the conditions of Theorems \ref{thm:main} or \ref{thm:tangency}.
\end{proposition}
\begin{proof}
Under the assumption on $H_1$ for $W$, we have $\omega$ is an integral class. To verify this, it suffices to check $\la \omega, B \ra\in \Z$ for $B\in H_2(X;\Z)$. We represent $B$ by a $f:\Sigma \to X$ from an oriented closed surface $\Sigma$. Then $f^{-1}(Y_{D})$ is a disjoint union of circles, whose image under the pushforward of $f$ represents a class in $H_1(Y_{D};\Z)$ that is trivial in $H_1(W;\Z)$. By our assumption, such a class is already zero in $H_1(Y_{D};\Z)$. In particular, the exactness of $\omega$ on $W$ implies that  $\la \omega, B \ra = \la \omega, f_*[\Sigma] \ra = \la \omega, g_*[\Sigma'] \ra$ for another closed oriented surface $\Sigma'$ and a map $g:\Sigma'\to N_{D}$ by replacing the null-homology in $W$ by the null-homology in $Y_{D}$. Now the claim follows from $g(\Sigma')\subset N_{D}\subset \CP^n$, where the symplectic form is the integer form $\omega_{\FS}$.  

Next, we can perturb $D_1$ to $D'_1$, such that $D'_1$ is a smooth hypersurface of the same degree. We may assume $D'_1$ is contained in $N_D$. Then the claim that $(X,\omega)$ satisfies the Gromov-Witten conditions in Theorems \ref{thm:main} and \ref{thm:tangency} follows from \Cref{lemma:GW} for $S=D'_1$.
\end{proof}

\begin{lemma}\label{lem:H1}
    Let $d=\gcd(\deg D_1,\ldots,\deg D_r)$. Then $H_1(Y_D;\Z)=\Z^{r-1}\oplus \Z/d$ generated by a small simple loop $\gamma_i$ around the smooth part of $D_i$ once, subject to the relation $\sum \deg D_i [\gamma_i]=0$.
\end{lemma}
\begin{proof}
    Let $\gamma$ be a loop in $Y_D$, we can find a bounding surface $v:\Sigma \to \CP^n$ with $v|_{\partial\Sigma}=\gamma$. By choosing $v$ generically such that it intersects the smooth parts of $D$, we obtain that $[\gamma]=\sum (v\cdot D_i) [\gamma_i]$ in $H_1(A_D;\Z)=H_1(Y_D;\Z)$. On the other hand, if $\sum a_i [\gamma_i]=0$, i.e.\ we have a surface $v:\Sigma \to A_D$, such that $v_*[\partial \Sigma]=\sum a_i[\gamma_i]$, by adding back the small disks intersect $D_i$ with multiplicity $a_i$, we obtain a map from a closed surface to $\CP^n$. Therefore we have $(a_1,\ldots,a_r)=k\cdot (\deg D_1,\ldots, \deg D_r)$ for $k\in \Z$. Hence the claim follows.
\end{proof}

We now move on to the proof of \Cref{thm:homology_unique}.
\begin{proof}[Proof of \Cref{thm:homology_unique}]
    In the standard decomposition $\CP^n=N_{D}\cup A_D$, where $N_D$ is the standard symplectic cap containing $D$, we have $A_D$ is a Weinstein domain. In particular, $\pi_{1}(Y_D)\to \pi_1(A_D)$ is an isomorphism. Hence, the van Kampen Theorem implies that $\pi_1(N_D)=1$. As a consequence, we have $\Tor H^2(N_D;\Z)$ is zero. 

    We write $(X=W\cup_{Y_D}N_D,\omega)$, \Cref{prop:integer} implies that $\omega$ is an integral class and $X$ is homotopy equivalent to $\CP^n$ by Theorems~\ref{thm:main} and \ref{thm:tangency}. We write $h\in H^2(X;\Z)$ as the unique integral lift of $\omega$, and define $h_N:=h|_{N_D}$. We write $h'_N$ as the restriction of the positive generator of $H^2(\CP^n;\Z)$ to $N_D$. Then we have $h'_N=h_N$. This is because, as real classes, they are both $\omega_{\FS}|_{N_D}$, and $H^2(N_D;\Z)$ is torsion-free. Now, since the cap $N_D$ contains a smooth complex hypersurface from smoothing all singularities, we have $h_N^j\ne 0$ rationally for $2j<2n$. That is $H^{j}(X:\Z)\to H^{j}(N_D;\Z)$ is injective for $j<2n$. Moreover, this map has the same image as $H^*(\CP^n;\Z)\to H^*(N_D;\Z)$. Now, note that 
    $$H_{k}(W;\Z)\simeq H^{2n-k}(W,Y_D;\Z)\simeq H^{2n-k}(X,N_D;\Z) \simeq \coker(H^{2n-k-1}(X;\Z)\to H^{2n-k-1}(N_D;\Z))$$
    which is $\coker(H^{2n-k-1}(\CP^n;\Z) \to H^{2n-k-1}(N_D;\Z))\simeq H_k(A_D;\Z)$. Here the third isomorphism follows from the injectivity of the said map and the long exact sequence for the pair. This confirms \ref{D1}. 

    For \ref{D2}, by Mayer-Vietoris and injectivity of $H^{j}(X;\Z)\to H^{j}(N_D;\Z)$ for $j<2n$, we have short sequence sequence: 
    $$0\to H^{j}(X;\Z)\to H^{j}(W;\Z)\oplus H^{j}(N_D;\Z)\to H^{j}(Y_D;\Z) \to 0,\quad j<2n-1.$$
    Hence $H^j(W;\Z)\to H^{j}(Y_D;\Z)$ is injective for $j<2n-1$. The injectivity for $j=2n-1$ is tautological by \ref{D1}. On the other hand, $H_{j}(W;\R)\to H_{j}(X;\R)$ is zero for $j>0$ by the pairing with $\omega^{j/2}$ and exactness of $W$. Since $H_*(X;\Z)$ is torsion-free, we have $H_{j}(W;\Z)\to H_{j}(X;\Z)$ is zero for $j>0$. Then the Mayer-Vietoris sequence
    $$\ldots \to H_j(Y_D;\Z)\to H_j(W;\Z)\oplus H_j(N_D;\Z)\to H_j(X;\Z) \to \ldots $$
    implies that every $(b,0)\in H_j(W;\Z)\oplus H_j(N_D;\Z),j>0$ is from $Y_D$, which proves the surjectivity. Finally as $A_D$ is Weinstein, we have $H_j(Y_D;\Z)\to H_{j}(A_D;\Z)$ is an isomorphism for $j\le n-2$ and is a surjection when $j=n-1$. The Mayer-Vietoris sequence implies that $H_j(N_D;\Z) \to H_j(\CP^n;\Z)$ is an isomorphism for $j\le n-2$ and surjection when $j=n-1$. We may apply the Mayer-Vietoris sequence to $X=W\cup_{Y_D}N_D$ to obtain $H_j(Y_D;\Z)\to H_{j}(W;\Z)$ is an isomorphism for $j\le n-2$. The isomorphism on cohomology then follows from the universal coefficient theorem.

    \ref{D3} follows from the van Kampen theorem as $\pi_1(N_D)=1,\pi_1(X)=0$.

    To see \ref{D4}, since $TN_D = T\CP^n|_{N_D}$ and the injection of $H^{j}(\CP^n;\Z)\to H^{j}(N_D;\Z)$ and $H^{j}(X;\Z)\to H^{j}(N_D;\Z)$  for $j<2n$, we have 
    $$c_j(TX) = \binom{n+1}{j} h^j.$$
    The top Chern class is also identified with $\CP^n$ as $\chi(X)=n+1$. Hence upon the identification between $H^*(X;\Z)$ and $H^*(\CP^n;\Z)$ as a ring, $X$ and $\CP^n$ have the same Chern character. Now as $H^*(X;\Z)$ is torsion free, \cite[Proposition 4.5]{KT} implies that $[TX]+[\underline{\C}] = (n+1)[L_h]$, where $L_h$ is the complex line bundle with first Chern class $h$. Since the restriction of $h$ to $Y_D$ is torsion of order $d$, and $\Tor H^2(W;\Z)\to \Tor H^2(Y_D;\Z)$ is an isomorphism. We have $h$ is a torsion class of order $d$ on $W$ as well.
\end{proof}

\subsection{Removing the $H_1$ condition}
In some cases, we can establish the injectivity on $H_1$ in \Cref{thm:homology_unique}. Recall that an irreducible algebraic hypersurface $D\subset \C^n$ is subcritical if $D=(f\circ \Phi)^{-1}(0)$, where $\Phi:\C^n\to \C^n$ is an polynomial isomorphism and $f$ is a non-trivial polynomial involving only $z_2,\ldots,z_n$. The significance of a subcritical hypersurface $D$ is that $\C^n\backslash D$ is subcritical. To see this, we may assume $\Phi=\Id$, then $\C^n\backslash D = \C_{z_1}\times (\C^{n-1}_{z_2,\ldots,z_n}\backslash f^{-1}(0))$, which is clearly subcritical. 

Throughout this part, $D_1,\ldots, D_r$ are $r$ subcritical hypersurfaces and $D=\cup D_i$ is reduced. By Lemma \ref{lem:H1}, we have $H_1(A_D;\Z)=\oplus^r \Z$ generated by small simple loops $\gamma_i$ around $D_i$. 

\begin{proposition}\label{prop:inj}
	Let $W$ be a strong filling of $Y=\partial (V \times \D)$ for a Weinstein domain $V$. Then $H_1(Y;\Q)\to H_1(W;\Q)$ is injective.
\end{proposition}
This was essentially established in \cite[Theorem 4.1 (b)]{BGMZ}. Although \cite[Theorem 4.1 (b)]{BGMZ} states for rational homology of degree $\dim_{\C}W-1$, the proof works for any degree $<\dim_{\C}W$. The degree constraint was for the purpose of applications in \cite{BGMZ} and was only used in  \cite[Theorem 4.1 (a)]{BGMZ} (but not optimal, see \cite[Remark 4.2]{BGMZ}). Briefly speaking, results in \cite{zbMATH07673358} imply that $\SH_{n-*}^{>0}(W;\Lambda)\to H^{*+1}(W;\Lambda)\to H^{*+1}(Y;\Lambda)$ is surjective onto the image of $H^*(DT^*S^1\times \D^k;\Lambda)\to H^*(Y;\Lambda)$ for the Novikov field $\Lambda$ over $\Q$ and any strong filling $W$. Here $\SH_{n-*}^{>0}(W;\Lambda)$ is the positive symplectic homology recalled in \S \ref{s:RFH}, but for strong fillings. Note that \cite{zbMATH07673358} used the symplectic cohomology convention. Hence, we have that $H^1(W;\Q)\to H^1(Y;\Q)$ is surjective, which is equivalent to $H_1(Y;\Q)\to H_1(W;\Q)$ is injective. 

\begin{proposition}\label{prop:H1}
    Let $W$ be a Liouville filling of $Y_D$, then $H_1(Y;\Z_D)\to H_1(W;\Z)$ is injective.
\end{proposition}
\begin{proof}
Assume otherwise, $0\ne \sum a_i[\gamma_i]\in H_1(Y_{D};\Z)$ is mapped to zero in $W$. W.O.L.G., we may assume $a_1\ne 0$, then we can fill back in hypersurfaces $D_{2},\ldots, D_{r}$ and glue with the Liouville filling $W$ of $Y_D$ to get a strong filling $V$ of $Y_{D_1}$. Since $Y_{D_1}$ is subcritically fillable, with $H_1(Y_{D_1};\Z)$ freely generated by $[\gamma_1]$. By \Cref{prop:inj},  $a_1[\gamma_1]$ is not null homologous in $V$. As a consequence,  $\sum a_i[\gamma_i]$ is not null homologous in $W$.
\end{proof}

\begin{example}[Complex hyperplane arrangements]\label{ex:HA}
Let $\cH = \{H_1,\dots,H_m\}$ be a finite collection of distinct hyperplanes in $\mathbb{C}^n$ for $n\ge 3$. For $D=\cup H_i$, $A_D$'s cohomology ring with integer coefficients is torsion‑free and is completely described by the combinatorial data of the arrangement. The classical result of Orlik and Solomon \cite{OS80} (for central arrangements), Orlik and Terao \cite[Theorem 5.90]{zbMATH00052944} (for general arrangements) gives an isomorphism of graded commutative algebras
\[
H^*(A_D;\Z) \;\cong\; OS(\cH),
\]
where $OS(\cH)$ is the \emph{Orlik–Solomon algebra} \cite[Definition 3.45]{zbMATH00052944} defined below.

Let $E = \bigwedge_{\mathbb{Z}} \langle e_1,\dots,e_m \rangle$ be the exterior algebra over $\mathbb{Z}$ on generators $e_i$ of degree $1$ corresponding to the hyperplanes $H_i$. For a subset $S = \{i_1,\dots,i_k\}$ with $i_1<\dots<i_k$, denote $e_S = e_{i_1}\wedge\cdots\wedge e_{i_k}$.  For each dependent subset $S\subseteq\{1,\dots,m\}$, i.e.\ $\bigcap_{i\in S} H_i\ne \emptyset$ has codimension $< |S|$, we define the boundary
\[
\partial(e_S) = \sum_{j=1}^{k} (-1)^{j-1} e_{S\setminus\{i_j\}}.
\]
Let $I$ be the ideal of $E$ generated by all $\partial(e_S)$ with $S$ dependent and $e_S$ for all $S$ with $\bigcap_{i\in S}H_i=\emptyset$. Then the Orlik–Solomon algebra is the quotient
\[
OS(\cH) = E/I.
\]
The algebra $OS(\cH)$ is graded by the exterior degree, and the isomorphism sends the class of $e_{i_1}\wedge\cdots\wedge e_{i_k}$ to the cohomology class of the logarithmic form $\frac{1}{2\pi \mathbf{i}}\frac{d\alpha_{i_1}}{\alpha_{i_1}}\wedge\cdots\wedge\frac{1}{2\pi \mathbf{i}}\frac{d\alpha_{i_k}}{\alpha_{i_k}}$, where $\alpha_i$ is  an affine linear function defining $H_i$. 
\end{example}

\begin{proof}[Proof of \Cref{thm:subcritical}]
    \ref{E1} follows from \Cref{prop:H1}. 

    For \ref{E2}, from \Cref{thm:homology_unique}, we only need to focus on degrees $n-1$ and $n$.  We need to look at two types of exact sequences from the Mayer-Vietoris sequence:
    $$0\to H^*(N_D;\Z)\oplus H^*(W;\Z)\to H^*(Y_D;\Z)\to 0, \quad *=n-1 \text{ or } n \text{ and is odd}$$
    $$0\to\Z \to H^*(N_D;\Z)\oplus H^*(W;\Z)\to H^*(Y_D;\Z)\to 0, \quad *=n-1 \text{ or } n \text{ and is even}$$
    By \Cref{thm:homology_unique} and \Cref{ex:HA}, $H^*(W;\Z)$ is torsion free, hence the map $\Z\to H^*(W;\Z)$ is zero. Therefore, the latter exact sequence is equivalent to
    $$0\to (H^*(N_D;\Z)/\Z)\oplus H^*(W;\Z)\to H^*(Y_D;\Z)\to 0, \quad *=n-1 \text{ or } n \text{ and is even}.$$
    By \Cref{thm:homology_unique}, $H^1(W;\Z)\to H^1(Y_D;\Z)$ is an isomorphism. Therefore the injective map $H^*(W;\Z) \to H^*(Y_D;\Z)$ shall cover at least the sub-algebra generated by $H^1(Y_D;\Z)$, which is isomorphic to $H^*(A_{D};\Z)$. Since this sub-algebra is already a complement to the image of $H^*(N_D;\Z)\to H^*(Y_D;\Z)$ and $H^*(N_D;\Z)/\Z\to H^*(Y_D;\Z)$ by looking at the standard case $W=A_D$. As a consequence, the two splittings above imply that the image of $H^*(W;\Z)\to H^*(Y_D;\Z)$ is the same as $\mathrm{im}(H^*(A_D;\Z)\to H^*(Y_D;\Z))$ for $*=n-1,n$.
\end{proof}

\subsection{Fillings of $ST^*\RP^n$}
In this part, we study Liouville fillings of $ST^*\RP^n$ equipped with the canonical contact structure for $n\ge 3$. It is the ideal contact boundary of $\CP^n\backslash Q$, where $Q$ is a smooth quadratic hypersurface in $\CP^n$. Let $N$ be a neighborhood of $Q$ as a symplectic cap and write $X=(W\cup_{Y}N,\omega)$ for a Liouville filling $W$ of $Y:=ST^*\RP^n$ with $n \ge 3$. Let $\gamma$ be a simple loop winding around $Q$ once. The first step towards \Cref{thm:RP} is proving the injectivity of $H_1(Y;\Z)\to H_1(W;\Z)$, i.e.\ $[\gamma]$ is non-trivial in $H_1(W;\Z)$.
\subsubsection{Injectivity on $H_1$ when $n\ge 4$}
Let $S$ be an embedded sphere in $Q$ with $\int_S \omega =1$. The homology class $[S]$ is the positive generator of $H_2(Q;\Z)=H_2(\CP^n;\Z)$. By Lemma \ref{lemma:GW}, we have 
$$\GW_{Q,1}(S,\pt) = \GW_{\CP^n,1}(S,\pt, 1) = 1.$$
By the same normal bundle trick in \cite[Lemma 2.3]{CieZhou}, we have 
$$\GW_{\CP^n,Q,1}(S,1,\pt) = \GW_{Q,1}(S,\pt)=1,$$
where $\GW_{\CP^n,Q,1}(S,1,\pt)$ counts rational curves in $\CP^n$ of degree $1$, passing through $S$ and a point, and is tangent to $Q$ (of order at least $1$) at the marked point passing through $S$. 

We require that the complex structure $J$ on $N$ is the standard one induced from $\CP^n$. 
We consider the following moduli spaces:
$$\cM_{S,\gamma}([S]):=\left\{u:\CP^1\to X\left| \begin{array}{l} \overline{\partial}_{J}u=0, u(0)\in \gamma, u(\infty)\in S, [u]=[S], \\ u \text{ tangent to  $Q$ at $\infty$} \end{array}\right.\right\}/_{\C^*}.$$
More precisely, by putting $\gamma$ near $Q$, we can assume curves in $\cM_{S,\gamma}([S])$ stays in $N$, as the limit case of shrinking $\gamma$ to a point on $Q$ would force all the curves to stay inside $Q$ by the positivity of intersection. In particular, $\cM_\gamma$ is an oriented compact $1$-dimensional manifold as the curve has minimal energy in $N$. $\GW_{\CP^n,Q,1}(S,1,\pt)=1$ implies that $\ev_0:\cM_{S,\gamma}([S])\to \gamma, u\mapsto u(0)$ is of degree $1$.

Now we assume that $[\gamma]=0$ in $H_1(W;\Z)$, we can find an oriented embedded surface $\Sigma\subset X$ disjoint from $Q$ and $\partial \Sigma=\gamma$. We define similarly 
\begin{equation}\label{eqn:moduli_sigma}
\cM_{S,\Sigma}([S]):=\left\{u:\CP^1\to X\left|\begin{array}{l} \overline{\partial}_{J}u=0, u(0)\in \Sigma, u(\infty)\in S, [u]=[S]\\ u \text{ tangent to  $Q$  at $\infty$} \end{array}\right.\right\}/_{\C^*}.
\end{equation}

\begin{claim}\label{claim:nobranch}
    For generic $J$, $\cM_{S,\Sigma}([S])$ contains no branched covered curves. In particular, $\cM_{S,\Sigma}([S])$ is a smooth oriented surface.  
\end{claim}
\begin{proof}
From the proof of Proposition \ref{prop:integer}, we see that $2[\omega]$ is integral. Hence, if a curve in $\cM_{S,\Sigma}([S])$ is branch covered, it is doubly curved. The underlying simple curves  shall have virtual dimension: 
        $$2n-2+2(n+1)/2-(2n-2)-(2n-2)=3-n < 0.$$
Note that any branch-covered curve must leave $N$, where we can perturb $J$ generically. Hence there are no branched covers for a generic $J$. 
\end{proof}
Let $\overline{\cM}_{S,\Sigma}([S])$ be the closure of  $\cM_{S,\Sigma}([S])$ inside the compactification by stable maps with constraint from $\Sigma $ and $S$ (but without tangency conditions). As any non constant $J$-curve in $X$ must interest $Q$ and the tangent point with $Q$ consumes up all the positive intersections of $[S]$ with $Q$, $\overline{\cM}_{S,\Sigma}([S]) \backslash \cM_{S,\Sigma}([S])$ only consists of nodal curves $(u_0,u_g,u_1)$ with 
\begin{enumerate}
    \item $u_0$ is a non-constant $J$-sphere with one marked point $0$, such that $u_0(0)\in S \subset Q$. We write $A=[u_0]\in H_2(X;\Z)$, then $A \cdot Q=1$ and $\omega(A) = \frac{1}{2}$.
    \item $u_1$ is a non-constant $J$-sphere with two marked points $0,\infty$, such that $u_1(0)=u_0(0)$ and $u_1(\infty)\in \Sigma$.  We write $B=[u_1]\in H_2(X;\Z)$, then $B \cdot Q=1$ and $\omega(B) = \frac{1}{2}$.
    \item $u_g$ is a ghost sphere at $u_0(0)$ with three marked points $(0,1,\infty)$, with $\infty$ connected with $0$ on $u_0$ and $0$ connected with $0$ on $u_1$.
\end{enumerate}
Note that $u_0,u_1$ have the minimal energy and leave $N$, we can compute the virtual dimensions of $u_0,u_1$ to conclude that for generic $J$, we must have 
$$(c_1(A),c_1(B)) = (1,n) \text{ or } (2,n-1).$$
And the nodal degeneration in $\overline{\cM}_{S,\Sigma}([S]) \backslash \cM_{S,\Sigma}([S])$ is of dimension $0$ and cut out transversely.

We first prove the following claim in several steps. 
\begin{claim}\label{claim:num0}
$\# (\overline{\cM}_{S,\Sigma}([S]) \backslash \cM_{S,\Sigma}([S]))=0$.
\end{claim}

Let $B$ be a class in $H_2(X;\Z)$ such that $c_1(B)=n$, $B\cdot Q=1$ and $\omega(B)=1/2$. We write $\overline{\Sigma}$ as $\Sigma$ glued with the small disk in $N$ bounding $\gamma$. We consider the moduli space 
$$\cM_{S,\overline{\Sigma}}(B):=\left\{u:\CP^1\to X\left|  \overline{\partial}_{J}u=0, u(0)\in \overline{\Sigma}, u(\infty)\in S, [u]=B \right.\right\}/_{\C^*}.$$
As those curves have minimal energy and leave $N$, for generic $J$, we can assume $\cM_{S,\overline{\Sigma}}(B)$ is a smooth oriented surface. Let $\cC\to \cM_{S,\overline{\Sigma}}(B)$ be the universal family, which is equipped with two sections $s_0,s_\infty$ from the marked points $0,\infty$. We may assume $\overline{\Sigma}\cap Q$ is disjoint from $S$. By positivity of intersection, we must have $s_0(\cM_{S,\overline{\Sigma}}(B))$ misses $\overline{\Sigma}\cap Q$. In particular, $(s_0)_*[\cM_{S,\overline{\Sigma}}(B)]$ is zero in $H_2(X;\Z)$. And we can assume $\cM_{S,\overline{\Sigma}}(B)=\cM_{S,\Sigma}(B)$.
\begin{claim}\label{claim:deg0}
    $s_\infty: \cM_{S,\overline{\Sigma}}(B) \to S$ has degree $0$.
\end{claim}
\begin{proof}
Assume the degree is $a\ne 0$, then $(s_\infty)_*[\cM_{S,\overline{\Sigma}}(B)]=a[S]$. Since $(s_0)_*[\cM_{\overline{\Sigma}}(B)]=0$, we have $(s_\infty)_*[\cM_{\overline{\Sigma}}(B)]=kB$ for $k\in \Z$ is the first Chern number of the line bundle behind $\cC$. From $kB=a[S]$ and the intersection with $Q$, we must have $k=2a\ne 0$. On the other hand, we have $2c_1(B)=2n \ne c_1(S)=n+1$. We arrive at a contradiction. 
\end{proof}
As a consequence, we have 
\begin{equation}\label{eqn:deg1}
    \# \left\{(u_0,u_g,u_1)\in \overline{\cM}_{S,\Sigma}([S]) \backslash \cM_{S,\Sigma}([S]), c_1(u_1)=n \right\}=0.
\end{equation}

Next, we consider the case with $c_1(u_1)=n-1$. Now let $B'$ be a class in $H_2(X;\Z)$ such that $c_1(B')=n-1$, $B'\cdot Q=1$ and $\omega(B')=1/2$. We write 
$$\cM_{S,\overline{\Sigma}}(B'):=\left\{u:\CP^1\to X\left|  \overline{\partial}_{J}u=0, u(0)\in \overline{\Sigma}, u(\infty)\in S, [u]=B' \right.\right\}/_{\C^*},$$
which is a compact oriented $0$-dimensional manifold for generic $J$. And we can also assume $\cM_{S,\overline{\Sigma}}(B')=\cM_{S,\Sigma}(B')$.

\begin{claim}\label{claim:count0}
    $\# \cM_{S,\overline{\Sigma}}(B')=0.$
\end{claim}
\begin{proof}
    Let $T=\CP^1\times \CP^1$ be a quadratic complex surface in $Q$, i.e.\ the intersection of $Q$ with a generic linear $\CP^3$. We may assume $S$ is contained in $T$ as $\CP^1\times \{\pt\}$. We consider 
    $$\cM_{T,\overline{\Sigma}}(B'):=\left\{u:\CP^1\to X\left|  \overline{\partial}_{J}u=0, u(0)\in \overline{\Sigma}, u(\infty)\in T, [u]=B' \right.\right\}/_{\C^*},$$
    which is an oriented closed surface for generic $J$. We write $C_1=[\CP^1\times \{\pt\}],C_2=[\{\pt\}\times \CP^1]\in H_2(T;\Z)$ and $(s_\infty)_*[\cM_{T,\overline{\Sigma}}(B')] = aC_1+bC_2$. Similarly to the proof of Claim \ref{claim:deg0}, as $C_1=C_2=[S]$ in $H_2(Q;\Z)$, we can conclude that $[S]=2B$ rationally if $a+b\ne 0$. Since $c_1(S)=n+1\ne 2c_1(B)=2(n-1)$, we must have $a+b=0$. As $C_1\cdot C_2=1$ on $T$, we have $\# \cM_{S,\overline{\Sigma}}(B')=b$. On the other hand, $\# \cM_{S,\overline{\Sigma}}(B')$ depends only on the homology class of $S$ in $X$. We may pick $S$ to be the $\{\pt\}\times \CP^1$ in $T$, as a consequence, we have $\# \cM_{S,\overline{\Sigma}}(B')=a$. Therefore, we must have $a=b=0$ and $\# \cM_{S,\overline{\Sigma}}(B')=0.$
\end{proof}

As a consequence, we have 
\begin{equation}\label{eqn:deg2}
    \# \left\{(u_0,u_g,u_1)\in \overline{\cM}_{S,\Sigma}([S]) \backslash \cM_{S,\Sigma}([S]), c_1(u_1)=n-1 \right\}=0.
\end{equation}

\begin{proof}[Proof of Claim \ref{claim:deg0}]
    It follows from the sum of \eqref{eqn:deg1} and \eqref{eqn:deg2}.
\end{proof}

In the following, we study the local picture near every curve in $\overline{\cM}_{S,\Sigma}([S]) \backslash \cM_{S,\Sigma}([S])$. We first study the neighborhood of $\overline{\cM}_{S,\Sigma}([S]) \backslash \cM_{S,\Sigma}([S])$ inside the compactified moduli spaces of stable maps without tangency conditions. 
Near $(u_0,u_g,u_1)\in \overline{\cM}_{S,\Sigma}([S]) \backslash \cM_{S,\Sigma}([S])$, we have a $\D_{\epsilon}\times \D_{\epsilon}$-chart into the compactified moduli space of stable maps. Let $z$ be the coordinate on $u_g$ and $z_0,z_1$ are coordinate on $u_0,u_1$. We write $t_0,t_1$ as the two gluing parameters from the two connected nodes from $u_0,u_1$, such that 
$$z_0z=t_0,\quad z_1/z=t_1.$$
Near $\mathrm{im}(u_g)\in Q$, we can write $Q$ locally as $\{w=0\}$ for a locally defined holomorphic function $w$, as $J$ is standard near $Q$. Therefore, we can write the glued curve on a central annulus of the middle component containing $z=1$ as the following Laurent series:
\begin{equation}\label{eqn:expansion}
    w_{\bm{t}}(z)=b_0(\bm{t})+\sum_{k>1}(b_k(\bm{t})z^k+b_{-k}(\bm{t})z^{-k})
\end{equation}
where $\bm{t}=(t_1,t_2)\in \D_{\epsilon}\times \D_{\epsilon}$ and $w_t(1)=0$. By our conditions, we have 
$$w\circ u_i(z_i)=\kappa_i z_i+O(z_i^2),\quad \kappa_i\ne 0,  i = 0,1.$$
If we plug $z=t_0/z_0$ into \eqref{eqn:expansion} and apply Cauchy integration formula to $\int_{|z_0|=1} w_t(t_0/z_0)/z_0^k \rd z_0$, we obtain that 
$$b_{-k}(\bm{t})/t_0^k \to \frac{1}{k!}(w\circ u_0)^{(k)}, t_0\to 0.$$
Similarly, we have 
$$b_{k}(\bm{t})/t_1^k \to \frac{1}{k!}(w\circ u_1)^{(k)}, t_1\to 0.$$
We define a map $$(p,q):\D_{\epsilon}\times \D_{\epsilon} \to \C^2, \quad \bm{t}\mapsto (w_t'(1), w_t(-1))$$
which has the asymptotic expansion 
$$p(\bm{t}) = -\kappa_0t_0+\kappa_1t_1+o(\bm{t}),\quad q(\bm{t})=w_t(-1)-w_t(1) =-2\kappa_0 t_0 -2\kappa_1 t_1+o(\bm{t}).$$
As a consequence, $\cM_{S,\Sigma}([S])$ near the nodal curve is described by $\kappa_0t_0=\kappa_1t_1$ approximately. 

\begin{remark}
In general, the compactification in relative symplectic Gromov-Witten theory is different from taking the closure inside the Gromov compactification of stable maps without tangency conditions. One needs to add in the rubber layers, see e.g.\ \cite[\S 2]{MZ}, \cite[\S 6]{IP} and \cite[\S 3]{LR}. In the special case here, the two compactification are identical with the ghost component viewed in the rubber layer.
The gluing Theorems in relative symplectic Gromov-Witten theory \cite{IP,LR,MZ} are harder as the gluing parameters from connecting nodes are no longer independent as above. The above straightforward argument endows the compactification $\overline{\cM}_{S,\Sigma}([S])$ a $C^0$ structure of surface at least, although what matters for us is just the linking number for a loop with $|\bm{t}|=\epsilon\ll 1$.
\end{remark}

Let $\eta$ be a small loop in $\cM_{S,\Sigma}([S])$ wrapping around a point in $\overline{\cM}_{S,\Sigma}([S])\backslash \cM_{S,\Sigma}([S])$ once. Let $C_{\eta}$ be the universal family over $\eta$ and $s$ be a section disjoint from $s_0$ of $C_{\eta}$. 
\begin{lemma}\label{lem:wind}
    For $T\gg 0$, the linking number of $Ts(\ell)$ with $Q$ is $1$ mod $2$.
\end{lemma}
\begin{proof}
    In the local model above, by choosing $s$ from the marked point $-1$ ($z=1$ is the marked point $\infty$ in \eqref{eqn:moduli_sigma}), we see that the linking number is precisely $1$. Note that changing $s$ will change the linking by a multiple of $2$, as $u$ is tangent to $Q$ with order precisely $1$ at $\infty$. Hence the claim holds in $\Z/2$.
\end{proof}

\begin{proposition}\label{prop:H_14}
    Let $W$ be a Liouville filling of $Y:=ST^*\RP^n$ for $n\ge 4$, then $H_1(Y;\Z)\to H_1(W;\Z)$ is an injection.
\end{proposition}
\begin{proof}
    Following the discussion above, let $\overline{\cM}$ be $\cM_{S,\Sigma}([S])$ deleting small disks around each point in $\overline{\cM}_{S,\Sigma}([S]) \backslash \cM_{S,\Sigma}([S])$. Then $\overline{\cM}$ is an oriented compact surface with boundary, we throw out all the closed components of $\overline{\cM}$ and still denote the remaining part by $\overline{\cM}$. Then the boundary of $\overline{\cM}$ are $\cM_{\gamma}$ and a small loop $\eta$ for each point in $\overline{\cM}_{S,\Sigma}([S]) \backslash \cM_{S,\Sigma}([S])$. Let $\overline{\cC}$ be the universal family over $\overline{\cM}$, which can be trivialized, as each component of $\overline{\cM}$ has boundary. Now let $s$ be a nowhere zero section of $\overline{\cC}$. We may consider the image of $Ts(\overline{\cM})$ for $T\gg 0$, which is contained in $N\backslash Q$. The small loop boundaries contribute to $0$ in $H_1(N\backslash Q;\Z/2)=\Z/2$ by Claim \ref{claim:num0} and Lemma \ref{lem:wind}. The boundary $\cM_{\gamma}$ contributes to $[\gamma]\ne 0 \in H_1(N\backslash Q;\Z/2)$ as all curves in $\cM_{\gamma}$ live in $N$ and $Ts$ is homotopic to $s_0$ without intersecting $s_\infty$. As a consequence, this shows that $[\gamma]=0\in H_1(N\backslash Q;\Z/2)$ by $Ts(\overline{\cM})$. We arrive at a contradiction. 
\end{proof}

\subsubsection{Injectivity on $H_1$ when $n=3$}
When $n=3$, both proofs of Claims \ref{claim:nobranch} and \ref{claim:count0} break down. Now $Q=\CP^1\times \CP^1$ and we write $C_1=[\CP^1\times \{\pt\}], C_2=[\{\pt\}\times \CP^1]$. We first assume $C_1\ne c_2\in H_2(X;\Z)$.Using the same normal bundle argument as in \cite[Lemma 2.3]{CieZhou}, we have 
$$\GW_{X,C_2}(C_1,C_1,\pt)=\GW_{Q,C_2}(C_1,C_1,\pt)=\GW_{X,C_1}(C_2,C_2,\pt)=\GW_{Q,C_1}(C_2,C_2,\pt)=1,$$
$$\GW_{X,C_1}(C_1,C_1,\pt)=\GW_{Q,C_1}(C_1,C_1,\pt)=\GW_{X,C_2}(C_2,C_2,\pt)=\GW_{Q,C_2}(C_2,C_2,\pt)=0,$$
$$\GW_{X,C_2}(C_1,C_2,\pt)=\GW_{Q,C_1}(C_1,C_2,\pt)=0.$$
For two classes $A,B\in H_2(X;\Z)$, we may represent them by disjoint surfaces in $X$. Then 
$\GW_{X,C_i}(A,B,\pt)\in \Z$, as there can not be any double branched covers in the moduli spaces by dimension reasons. As a consequence, we have $C_i$ is primitive in $H_2(X;\Z)$ modulo torsion. For otherwise, if $C_i=kA$ modulo torsion in $H_2(X;\Z)$, we have 
$$k^2\GW_{X,C_{1-i}}(A,A,\pt)=1$$
and hence $k=\pm 1$. Similar to \eqref{eqn:moduli_sigma}, we can define $\cM_{S,\gamma}([C])$ and $\cM_{S,\Sigma}([C])$ using curves in class $C=C_i$ and $S\subset Q$ representing $C_{1-i}$. The above discussion implies that $\cM_{S,\Sigma}([C])$ is an oriented smooth surface. To study the compactification, the analogue of Claim \ref{claim:deg0} holds. To see analogue of Claim \ref{claim:count0}, by the same argument we have 
$$aC_1+bC_2=2(a+b)B \in H_2(X;\Z)$$
But we can not argue that $a=b$, we do not know $c_1=C_2$ in $H_2(X;\Z)$. Note that we have 
$$b^2=\GW_{X,C_1}(aC_1+bC_2,aC_1+bC_2,\pt) =4(a+b)^2 \GW_{X,C_1}(B,B,\pt)\in 4(a+b)^2\Z, $$
$$a^2=\GW_{X,C_2}(aC_1+bC_2,aC_1+bC_2,\pt) =4(a+b)^2 \GW_{X,C_2}(B,B,\pt)\in 4(a+b)^2\Z $$
If $a=0$, then we must have $b=0$ from the first equation. Similarly if $b=0$, we must have $a=0$. And if $a+b=0$, then $a=b=0$. Next we assume $a,b,a+b\ne 0$. Then we have $2$ divides both $a,b$, and inductively $2^k$ divides $a,b$ for any $k\in \N$, hence no such $a,b$ exist. As a consequence, we can conclude that $a=b=0$. Hence we have the analogue of Claim \ref{claim:count0}. 

Now if $C_1=C_2=C$ in $H_2(X;\Z)$, then we have 
$$\GW_{X,C}(C,C,\pt) = \GW_{Q,C_2}(C_1,C_1,\pt)=1.$$
The same argument implies that $C$ not not be decomposed integrally. For the analogue of Claim \ref{claim:count0}, we can argue as in the $n\ge 4$ case.

As a consequence, we have the following by the same argument using Lemma \ref{lem:wind}. 
\begin{proposition}\label{prop:H13}
Let $W$ be a Liouville filling of $Y:=ST^*\RP^3$, then $H_1(Y;\Z)\to H_1(W;\Z)$ is an injection.   
\end{proposition}

\subsubsection{The proof of \Cref{thm:RP}}

\begin{proposition}\label{prop:ring}
    Let $D$ be a smooth complex hypersurface of degree $k<n$ in $\CP^n$. Let $W$ be a Liouville filling of $Y_D$, such that $H_1(Y_D;\Z)\to H_1(W;\Z)$ is injective, then $H^*(A_D;\Z)$ is isomorphic to $H^*(W;\Z)$ as a graded ring.
\end{proposition}
\begin{proof}
By Lefschetz hyperplane theorem, $H^*(\CP^n;\Z)\to H^*(D;\Z)$ is an isomorphism for $*<n-1$ and an injection for $*=n-1$. Moreover, $H^*(D;\Z)$ is torsion free. The Thom isomorphism and the exact sequence of the pair $(\CP^n, A_D)$ imply that 
\begin{equation}\label{eqn:Thom}
    \ldots H^{j-2}(D;\Z)\to H^{j}(\CP^n;\Z)\to H^{j}(A_D;\Z) \to H^{j-1}(D)\to H^{j+1}(\CP^n;\Z)\ldots
\end{equation}
Such that $H^{j-2}(D;\Z)\to H^j(\CP;\Z)$ is the cup product with $kh$ for the positive generator $h\in H^2(\CP^n;\Z)$ when $j\le n+1$. Therefore, from \eqref{eqn:Thom}, we can deduce that $H^*(A_D;\Z)=\Z[h]/\la kh, h^{\lfloor n/2\rfloor+1} \ra\oplus \Z^r[-n]$ as a graded ring, where $r=b_{n-1}(D):=\rank H^{n-1}(D;\Z)$ when $n$ is even and $r=b_{n-1}(D)-1$ when $n$ is odd. Since \eqref{eqn:Thom} also holds for $X=W\cup_{Y_D}N_D\simeq \CP^n$ and $H^{j-2}(D;\Z)\to H^{j}(X;\Z)$ is also identified with $H^{j-2}(D;\Z)\to H^{j}(\CP^n;\Z)$ by the proof of \Cref{thm:homology_unique}, hence we can deduce that $H^*(W;\Z)$ is also $\Z[h]/\la kh, h^{\lfloor n/2\rfloor+1} \ra\oplus \Z^r[-n]$ as a graded ring.
\end{proof}

The proof of \ref{F2} of \Cref{thm:RP} requires the following result by Kwon-Oba \cite{zbMATH08152684} and  Bimmermann-Stratmann-Zehmisch \cite{arXiv:2506.07525}.

\begin{theorem}[Kwon-Oba, Bimmermann-Stratmann-Zehmisch]\label{thm:Sn}
Let $W$ be a Liouville filling of $ST^*S^n$ for $n\ge 3$ odd. We consider a section $s$ of $T^*S^n$ with $|s|=1-\epsilon$. Then $W\backslash N(s)$ is an $h$-cobordism from $ST^*S^n$ to itself, where $N(s)\simeq DT^*S^n$ is a sufficiently small tubular neighborhood of the section $s$. In particular, $W$ is diffeomorphic to $DT^*S^n$ by the $h$-cobordism theorem. 
\end{theorem}

To evoke \Cref{thm:Sn}, we need to obtain that $\pi_1(ST^*\RP^n)\to \pi_1(W)$ is an isomorphism. This does not follow from \Cref{thm:homology_unique}, Propositions \ref{prop:H_14} and \ref{prop:H13}. Indeed, we could have 
$$\Z/2 \stackrel{\iota}{\to} \pi_1(W)\stackrel{\mu}{\to} \pi_1(W)/[\pi_1(W),\pi_1(W)]=\Z/2$$
with $\mu\circ \iota$ is an isomorphism and $\iota$ normally generates $\pi_1(W)$, e.g.\ we could have $\pi_1(W)$ is the permutation group $S_3$ and $\iota$ sends the non-trivial element to the transposition $(12)$.

We first prove the following more general result regarding the fundamental group.
\begin{proposition}\label{prop:pi1}
    Under the assumption of \Cref{thm:homology_unique} with $D$ a smooth complex hypersurface of degree $k<n$, we have $\pi_1(Y_D)\to \pi_1(W)$ is an isomorphism.
\end{proposition}
\begin{proof}
    By \Cref{thm:homology_unique}, we have $\pi_1(Y_D)\to H_1(Y_D;\Z)\simeq H^1(W;\Z)=\Z/k$ is an isomorphism. We can consider the $k$-fold covering $\widetilde{W}$ of $W$ from the normal subgroup $[\pi_1(W),\pi_1(W)]$. Then $\widetilde{W}$ is a Liouville filling of the universal cover $\widetilde{Y}_D$, which is the degree $1$ circle bundle over a smooth degree $k$ complex hypersurface $D$ in $\CP^n$ equipped with the Boothby-Wang contact structure. We claim that $\pi_1(\widetilde{W})=1$, which is equivalent to that $[\pi_1(W),\pi_1(W)]=0$ and hence $\pi_1(Y_D)\to \pi_1(W)$ is an isomorphism. Let $M$ be $\widetilde{W}$ capped by adding by the divisor $D$, which is a closed manifold symplectic manifold and $\pi_1(M)=\pi_1(\widetilde{W})$. Now we consider the universal cover $\widetilde{M}$ of $M$. If $\pi_1(M)\ne 1$, then $\widetilde{M}$ contains multiple copies of $D$ with normal bundle of degree $1$. Let $D_1,D_2$ be the two copies. By \cite[Lemma 2.3]{CieZhou} and Lemma \ref{lemma:GW}, we get 
    $$0\ne  \sum_A \GW_{D_1,A}(B,\pt) = \GW_{\widetilde{M},A}(B,\pt)$$
    where $A$ is any homology class from $D_1$ corresponding to the positive generator in $\CP^n$\footnote{There is exactly one such class unless $(n,k)=(3,2)$, where there are two such classes with nonzero $\GW$.} and $B$ is $2(k-1)$ dimensional submanifold of $D_1$ such that $\la B,\omega^{k-1} \ra\ne 0$. Although $\widetilde{M}$ might be non-compact, if we use an almost complex structure on $\widetilde{M}$ pulled back from $M$, we still have Gromov compactness to define the Gromov-Witten invariants for passing through $B$ and a point by the unique lifting property of the covering. Now, if we put the point constraint in $\GW_{\widetilde{M},A}(B,\pt)$ on $D_2$, since $A\cdot D_2=0$, the positivity of intersection implies that the curve for $\GW_{\widetilde{M},A}(B,\pt)$ must be contained in $D_2$, contradicting with $B\subset S_1$. This concludes the proof.
\end{proof}

\begin{proof}[Proof of \ref{F1} and \ref{F2} of \Cref{thm:RP}]
\ref{F1} follows from \Cref{thm:homology_unique}, Propositions \ref{prop:H_14}, \ref{prop:H13}, \ref{prop:pi1} and \ref{prop:ring}. The last claim follows from the proof of \ref{D4} in \Cref{thm:homology_unique} and the fact that the generator of $H^2(\CP^2;\Z)$ restricted to the generator of $H^2(ST^*\RP^n;\Z)=\Z/2$ which is isomorphic to $H^2(W;\Z)$ by \Cref{thm:homology_unique}.

For \eqref{F2}, when $n$ is odd, we can find a nowhere-zero section $s$ of $DT^*\RP^n$ with $|s|=1-\epsilon$, hence we can move a neighborhood of the zero section into the neighborhood of $ST^*\RP^n$. We have a cobordism $X=W\backslash N(s)$ between $ST^*\RP^n$ for a sufficiently small tubular neighborhood $N(s)$ of the section $s$. By \Cref{prop:pi1}, $\pi_1(ST^*\RP^n)\to \pi_1(W)$ is an isomorphism. Then the universal cover $\widetilde{W}$ is a Liouville filling of $ST^*S^n$. Then the van Kampen Theorem implies that $\pi_1(X)=\Z/2$ and the inclusion of each of the boundary components induces an isomorphism on the fundamental group. Moreover, the universal cover $\widetilde{X}$ of $X$ is precisely the $h$-cobordism in \Cref{thm:Sn}. Therefore, $X$ is also an $h$-cobordism. Since the Whitehead group of $\Z/2$ is trivial, $s$-cobordism theorem implies that $X$ is a trivial cobordism. Therefore, $W$ is diffeomorphic to $DT^*\RP^n$.
\end{proof}

For \ref{F3} of \Cref{thm:RP}, we divide the proof into several steps. Let $W$ be a Weinstein filling of $Y:=ST^*\RP^n$ for $n\ge 4$ even. The idea is to first build a homotopy equivalence $\RP^n\to W$ along with a tangential map from \ref{F1} of \Cref{thm:RP}, then use the Smale-Hirsch immersion Theorem \cite{Hirsch} to homotope the map through formal immersions into an actual immersion, finally use the Hatcher-Quinn Theorem \cite{HQ} to update it to an embedding through immersions. This allows us to embed $DT^*\RP^n$ into $W$ in a homotopy equivalent way, and we finish the proof by the usual $h$-cobordism argument. 

\begin{claim}\label{claim:double}
    The universal cover $\widetilde{W}$ is homotopy equivalent to $S^n$.
\end{claim}
\begin{proof}
    The universal cover also has a skeleton of dimension $n$, as $\partial \widetilde{W} = ST^*S^{n}$ is $n-2$ connected, we have $\widetilde{W}$ is also $n-2$-connected. From the exact sequence
    $$\Z/2=H^{n}(ST^*S^n;\Z)\to H^{n+1}(\widetilde{W},\partial \widetilde{W};\Z) \to H^{n+1}(\widetilde{W};\Z)=0$$
    and Lefschetz duality $H_{n-1}(\widetilde{W};\Z) = H^{n+1}(\widetilde{W},\partial{\widetilde{W}};\Z)$, we have $H_{n-1}(\widetilde{W};\Z)$ is either $\Z/2$ or $0$. Since $\chi(\widetilde{M})=2\chi(W)=2$, we have $H^n(\widetilde{W};\Z)$ has rank one. If $H_{n-1}(\widetilde{W};\Z)=\Z/2$, by universal coefficient theorem, we have $H^{n}(\widetilde{W};\Z)=\Z\oplus \Z/2$. Hence we have exact sequence 
    $$0\to H^n(\widetilde{W},\partial \widetilde{W};\Z)\to H^n(\widetilde{W};\Z)\to H^n(\partial\widetilde{W};\Z)\to H^{n+1}(\widetilde{W},\partial\widetilde{W};\Z)\to 0$$
    which is 
    $$0\to \Z \to \Z\oplus \Z/2 \to \Z/2 \to \Z/2\to 0.$$
    We arrive at a contradiction. Therefore $\widetilde{W}$ is simply connected and has the same homology as $S^n$, hence is homotopy equivalent to $S^n$. 
\end{proof}
Claim \ref{claim:double} is the only place where the Weinstein property is used. 

\begin{claim}\label{claim:HE}
    There exists a homotopy equivalence $f:\RP^n \to W$, such that $f^*TW \simeq TP\oplus T^*P$.
\end{claim}
\begin{proof}
    Let $\gamma$ be a loop in $W$ representing the generator of $\pi_1(W)$. Viewing $\gamma$ as the $1$-Skelton of $\RP^n$, since $\gamma^2$ is contractible, we can extend the map to $\RP^2 \to W$. By Claim \ref{claim:double}, $\pi_i(W)=0$ for $2\le i <n$, we can extend the map to $f:\RP^n\to W$. Since $f$ induces isomorphism on $H_1(\cdot;\Z)$, $f$ induces an isomorphism on $H^2(\cdot;\Z)$. By \ref{F1} of \Cref{thm:RP}, $f$ induces an isomorphism on both $H^*$ and $H_*$ as $H^*(W;\Z)$ is generated by $H^2(W;\Z)$. Hence $f$ induces isomorphism on (co)homology for any coefficient. Now let $\widetilde{f}:S^n\to \widetilde{W}\simeq S^n$ be the lift of $f$. If we consider the spectral sequence computing $H^*(W;\Z/2)$ by the fiberation $\widetilde{W} \to \widetilde{W}\times_{\Z/2}E\Z/2 \to B\Z/2$, we have a non-trivial differential from $H^n(\widetilde{W};\Z/2)=\Z/2 \to H^{n+1}(B\Z/2;\Z/2)=\Z/2$. Now the compatibility of the spectral sequence for $\Z/2 \circlearrowright S^n$ using $\widetilde{f}$, implies that $\deg \widetilde{f}$ is odd. Therefore, we can change the map on the last cell of $\RP^n$ to make $\widetilde{f}$ have degree exactly $1$. In particular, this $f$ is a homotopy equivalence. Now by \ref{F1} of \Cref{thm:RP},$[f^*TW]\oplus [\underline{\C}]=(n+1)[L_h]$, where $h$ is the generator on $H^2(\RP^n;\Z)$. Hence as real bundles, we have $[f^*TW]=[TP\oplus T^*P]$ stably. Now by \cite[Chapter 9, Theorem 1.5]{fiber}, we can get an honest bundle isomorphism $f^*TW \simeq TP\oplus T^*P$. 
\end{proof}

\begin{claim}\label{claim:embedding}
    There exists an embedding $\iota:\RP^n \to W$, such that $\iota$ is a homotopy equivalence and the normal bundle is isomorphic to $T^*P$. 
\end{claim}
\begin{proof}
By Claim \ref{claim:HE}, we have a formal immersion $(f,F)$, where $F$ is the inclusion of $T\RP^n$ into $f^*TW$ via the isomorphism to $T\RP^n\oplus T^*\RP^n$. By the Smale-Hirsch immersion Theorem \cite{Hirsch}, we can homotope $(f,F)$ in the space of formal immersions to an actual immersion $g:\RP^n \to W$ with only transverse double points. Let $\mathfrak{o}$ be the non-trivial sign local system on $W$, then $g^*\mathfrak{o}$ is the determinant bundle of $T\RP^n$. Since $H_n(\RP^n;g^*\mathfrak{o})=H^n(\RP^n;g^*\mathfrak{o})=\Z$, we have $H_n(W;\mathfrak{o})=H^n(W;\mathfrak{o})=\Z$. And $g$ sends the generator of $H_n(\RP^n;g^*\mathfrak{o})$ to a generator $x$ of $H_n(W;\mathfrak{o})$. Since $H_{n-1}(ST^*\RP^n;\mathfrak{o})=0$, e.g.\ from the Gysin sequence computation, we have $H_n(W;\mathfrak{o})\stackrel{i}{\to} H_n(W,\partial W;\mathfrak{o})$ is surjective and hence an isomorphism. As a consequence, we have the intersection form on $H_n(W;\mathfrak{o})$, i.e. $H_n(W;\mathfrak{o})\otimes H_n(W;\mathfrak{o})\to \Z$ defined by $\la i(\cdot),\cdot \ra$, is unimodular. In particular, $x\cdot x = \pm 1$. Now for an immersion $g:\RP^n \to W$, we obtain an exact sequence $0\to T\RP^n\to g^*TW \to N\to 0$, where $N$ is the normal bundle. Since $W$ is oriented, we get an Euler number $N(g)\in \Z$ of the normal bundle. In our case here, we have $N(g)=\pm \chi(\RP^n)=\pm 1$. In the construction of $F$ in the formal immersion, we may compose an automorphism of $f^*TW$ that changes the orientation, this allows us to change the sign of $N(g)$ for the immersion produced by the Smale-Hirsch immersion Theorem. As a consequence, we may assume $N(g)=x\cdot x$.

Finally, we need to regularly homotope $g$ to an embedding. This requires using the Hatcher-Quinn Theorem \cite{HQ}, we refer readers to \cite[\S 2,3]{Grant} for a simpler account.  Let $\cP:=\RP^n \times^h_W \RP^n$ be the homotopy pullback, namely tuples $(x,\eta,y)$ for $x,y\in \RP^n$ and $\eta$ is a path in $W$ connecting $g(x)$ and $g(y)$. Since $h$ is a homotopy equivalence, the natural diagonal $\RP^n \to \cP$ is also a homotopy equivalence, and is equivariant w.r.t.\ the $\Z/2$ action in \cite[(25)]{Grant} on $\cP$ and trivial action on $\RP^n$. Let $p_1,p_2$ be two projections to $N$ and $p:\cP\to W$ sends $(x,\eta,y)$ to $\eta(1/2)$. Let $\zeta$ denotes the virtual bundle $p_1^*T\RP^n\oplus p_2^*\RP^n - p^*TW$, which also carries a $\Z/2$ action. $\zeta$ over the diagonal $\RP^n$ is orientable, as $w_1(\zeta|_{\RP^n})=2w_1(\RP^n)-w_1(g^*TW)=0$.  Since $\Z/2$ action on $\xi$ over the diagonal preserves orientation, we have $\xi\times_{\Z/2} {E\Z/2}$ is orientable over $P\times B\Z/2$.  By \cite[(6) and Theorem 3.3]{Grant}, the obstruction to embedding $\mu(g)$ lives in $\Omega_0(\cP\times_{\Z/2}E\Z/2, \zeta\times_{\Z/2}E\Z/2)=H_0(\RP^n\times B\Z/2;\Z)$ via the diagonal inclusion. Tracing the definitions in \cite[\S3]{Grant}, $\mu(g)$ is precisely the count of unordered double points. Here the sign is defined as follows: if $g(x)=g(y),x\ne y$, choosing an orientation in $\mathfrak{o}(g(x))$ on $W$ determines an orientation of $T_x\RP^n$ and $T_y\RP^n$, the intersection is positive if the orientations sum to the orientation of $T_{g(x)}W$. Now it is easy to check that $x\cdot x = E(g)+2\mu(g)$, our construction yields that $\mu(g)=0$. Hence $g$ is regularly homotopic to an embedding $\iota$ with normal bundle $T^*\RP^n$.
\end{proof}

\begin{proof}[Proof of \ref{F3} of \Cref{thm:RP}]
    By Claim \ref{claim:embedding}, we can remove the neighborhood $N(\iota)$ of the embedding to obtain a cobordism $X$ from $ST^*\RP^n$ to itself. $\iota$ being a homotopy equivalence implies that both boundary inclusions of $X$ induce isomorphisms on the fundamental group. Now we can consider the double covers of everything. $\iota$ being a homotopy equivalence implies that $\widetilde{N(\iota)}\hookrightarrow \widetilde{W}$ are homotopy equivalences for the double cover. Hence $H_k(\widetilde{X},\partial \widetilde{N(\iota)};\Z)=H_k(\widetilde{W},\widetilde{N(\iota)};\Z)=0$. This proves that $\partial N(\iota)\to X$ is a homotopy equivalence. Using duality and the universal coefficient theorem, we can deduce that the other boundary inclusion is also a homotopy equivalence. Hence $X$ is an $h$-cobordism and an $s$-cobordism as $\pi_1(X)=\Z/2$. Therefore $W$ is diffeomorphic to $DT^*\RP^n$.
\end{proof}

\bibliographystyle{alpha} 
\bibliography{ref}
\Addresses

\end{document}